%% file: main.tex
\documentclass[conf ]{new-aiaa}
\usepackage[utf8]{inputenc}
\usepackage{textcomp}

\usepackage{graphicx}
\usepackage{color}
\usepackage{amsmath, amsfonts}
\usepackage{amssymb}
\usepackage[version=4]{mhchem}
\usepackage{siunitx}
\usepackage{longtable,tabularx}

\usepackage{bm}
\usepackage{subfigure}
\usepackage{wasysym}
\usepackage{algorithm}
\usepackage{multirow}

\usepackage{braket}
\usepackage{doi}
\usepackage{url}
\usepackage{fancyhdr}

\title{Tube Stochastic Optimal Control for Nonlinear Constrained Trajectory Optimization Problems}

\author{Naoya Ozaki\footnote{Ph.D., Assistant Professor, Department of Spacecraft Engineering, Japan Aerospace Exploration Agency, Kanagawa 252-5210, Japan; ozaki.naoya@jaxa.jp, and Member AIAA.}}
\affil{Japan Aerospace Exploration Agency, Sagamihara, Kanagawa, 252-5210, Japan}

\author{Stefano Campagnola\footnote{Mission Design Engineer, Outer Planet Mission Analysis Group, Jet Propulsion Laboratory, California Institute of Technology, Pasadena, CA, 91109, USA.}}
\affil{Jet Propulsion Laboratory, California Institute of Technology, Pasadena, CA, 91109, USA}

\author{Ryu Funase\footnote{Ph.D., Associate Professor, Department of Aeronautics and Astronautics, the University of Tokyo, 7-3-1, Hongo, Bunkyo-ku, Tokyo, 113-8656, Japan, and Member AIAA.}}
\affil{The University of Tokyo, Bunkyo-ku, Tokyo, 113-8656, Japan}

\begin{document}

\maketitle

\input{contents/0_abstract}

\section*{Nomenclature}

{\renewcommand\arraystretch{1.0}
\noindent\begin{longtable*}{@{}l @{\quad=\quad} l@{}}
$\mathbb{E}$ & expected value\\
$k$ & stage (node) number; $k\in\mathbb{N}$\\
$l$ & stage cost functions; $l: \mathbb{R}^{n_x}\times \mathbb{R}^{n_u}\times \mathbb{R}^{n_w} \rightarrow \mathbb{R}$\\
$\mathcal{N}(\bm{m},\bm{\mathrm{S}})$ & Gaussian distribution with mean $\bm{m}$ and covariance $\bm{\mathrm{S}}$\\
$\mathbb{N}_{b}$ & set of non-negative integers from $0$ to $b$; $\mathbb{N}_b:=\{0,1,...,b\}$\\
$\mathbb{N}_{a:b}$ & set of integers between $a$ and $b$; $\mathbb{N}_{a:b}:=\{a,a+1,...,b\}$\\
$N$ & number of stages; $N\in\mathbb{N}$\\
$n_x$ & dimension of state; $n_x\in\mathbb{N}$\\
$n_u$ & dimension of control; $n_u\in\mathbb{N}$\\
$n_w$ & dimension of uncertainty; $n_w\in\mathbb{N}$\\
$\bm{\mathrm{P}}$ & covariance matrices of state vector; $\bm{\mathrm{P}}\in\mathbb{R}^{n_x\times n_x}$\\
$\bm{\mathrm{R}}$ & covariance matrices of uncertainty; $\bm{\mathrm{R}}\in\mathbb{R}^{n_w\times n_w}$\\
$\mathbb{V}$ & variance \\
$\bar{\bm{x}}$ & mean values of state vector; $\bar{\bm{x}}\in\mathbb{R}^{n_x}$\\ 
$\varphi$ & terminal cost function; $\varphi: \mathbb{R}^{n_x} \rightarrow \mathbb{R}$\\
$\bm{0}_n$ & $n$-dimensional zero vector\\
$O_n$ & $n\times n$-dimensional zero matrix
\end{longtable*}}

%
%
\input{contents/1_introduction}

\input{contents/2_socp}

\input{contents/3_tsoc_by_ut}
\input{contents/4_num_example_double_integ}
\input{contents/5_num_example_low_thrust}
\input{contents/6_conclusion}
\input{contents/9_appendix}

%
%
\section*{Acknowledgments}\label{Acknowledgments}

This work was supported by the Japan Society for the Promotion of Science (JSPS) Grant-in-Aid for JSPS Fellows Number 15J05999/18J02128 and JPL Visiting Student Researchers Program (JVSRP). Part of this research was carried out at the Jet Propulsion Laboratory, California Institute of Technology, under a contract with the National Aeronautics and Space Administration. The first author would like to thank Jon Sims [Jet Propulsion Laboratory,California Institute of Technology(JPL)], Gregory Lantoine, Gregory Whiffen (JPL), Jonathan Aziz (CU Boulder), Nathaniel Guy (JPL), Eric Gustafson (JPL), Nicola Baresi (JAXA), Ferran Gonzalez-Franquesa, and Onur Celik (SOKENDAI, the Graduate University for Advanced Studies) for their valuable comments. This work was presented as Paper 2019 at the 29th AAS/AIAA Space Flight Mechanics Meeting in Maui, Hawaii, 13-17 January 2019.

\bibliography{references}

\end{document}

%% file: contents/0_abstract.tex
\begin{abstract}
Recent low-thrust space missions have highlighted the importance of designing trajectories that are robust against uncertainties. In its complete form, this process is formulated as a nonlinear constrained stochastic optimal control problem. This problem is among the most complex in control theory, and no practically applicable method to low-thrust trajectory optimization problems has been proposed to date. This paper presents a new algorithm to solve stochastic optimal control problems with nonlinear systems and constraints. The proposed algorithm uses the unscented transform to convert a stochastic optimal control problem into a deterministic problem, which is then solved by trajectory optimization methods such as differential dynamic programming. Two numerical examples, one of which applies the proposed method to low-thrust trajectory design, illustrate that it automatically introduces margins that improve robustness. Finally, Monte Carlo simulations are used to evaluate the robustness and optimality of the solution.
\end{abstract}

%% file: contents/1_introduction.tex
%
%
\section{Introduction}

%
%
%
Stochastic optimal control, which controls systems with probabilistic uncertainties, is of interest from both theoretical and practical perspectives. Over the past three decades, many researchers have studied methods to solve stochastic optimal control problems. However, numerical algorithms to solve nonlinear constrained problems, in which trajectory optimization problems are formulated, are an area of research that remains largely unexplored. 

%
%
%
In constrained systems, the presence of uncertainties can drive the optimal control problem to infeasibility. Chisci et al.\cite{Chisci2001} introduced the idea of improving robustness against disturbances by enforcing suitable constraint restrictions to the control inputs, referred to as the \textit{constraint tightening} method\cite{Chisci2001, Richards2006}. This is analogous to introducing a \textit{duty cycle} into a low-thrust trajectory design\cite{Rayman2007}. Mayne et al.\cite{Mayne2005} successfully implemented the idea of \textit{tube model predictive control}\cite{Langson2004, Cannon2009, Rakovic2012, Rakovic2012a} for linear time-invariant systems by introducing a tube created by a \textit{robust positively invariant set}\cite{Rakovic2007}. While their method is one of the most widely used methods to solve constrained problems, it is usually limited to linear quadratic problems.


In nonlinear problems, stochastic optimal control problems are generally formulated by the \textit{Bellman equation}\cite{Bertsekas1996, Bertsekas2000} or the \textit{Hamilton-Jacobi-Bellman equation} (HJB equation)\cite{Yong1999, Bertsekas2000}. Some methods solve the partial differential equations of the HJB equation directly using numerical techniques\cite{Yong1999, Kappen2005, Gustafson2010}; however, the methods incur a substantial computational cost because of the high-dimensional search space. Other methods iteratively solve the locally expanded Bellman equation around a reference trajectory\cite{Todorov2005, Todorov2009, Theodorou2010, Ozaki2018}. 
Notable local expansion algorithms include \textit{Stochastic Differential Dynamic Programming} (SDDP)\cite{Theodorou2010, Boutselis2016, Ozaki2018}, which is based on \textit{Differential Dynamic Programming} (DDP)\cite{Jacobson1970, Lantoine2012p1}. DDP solves a second-order expansion of the Bellman equation to find the local optimal control. Ross et al.\cite{Ross2014, Ross2015} have proposed another approach, known as the \textit{unscented optimal control} algorithm. This algorithm introduces the \textit{unscented transform}\cite{Julier1996, Julier1997} to model uncertainties in a deterministic way so that the resultant problem can be solved by a deterministic method\cite{Ross2007}. Nonetheless, these approaches are limited to open-loop control without correcting perturbations, and cannot produce results that are as robust as the constraint tightening method.\cite{Chisci2001}

%
%
%
%
Currently, none of these approaches have resulted in an algorithm that can model stochastic systems in a nonlinear way, handle the closed-loop control with constraints, and solve high-dimensional optimization problems. 

%
%
%
This paper presents a novel algorithm to solve stochastic optimal control problems that possess nonlinearities and constrain control inputs in a probabilistic manner. The proposed method sequentially approximates stochastic processes as Gaussian processes, employs the unscented transform to estimate the expected cost and stochastic evolution of the dynamics, and applies a chance-constrained method to handle the control constraints. These approximations transcribe the stochastic optimal control problem into a deterministic problem. The deterministic problem is then solved by trajectory optimization methods such as DDP.

%% file: contents/2_socp.tex
%
%
\section{Stochastic Optimal Control Problem}

This section summarizes the basics of stochastic optimal control theory, and introduces the notation and assumptions required to derive the proposed method. Our model conforms to a discrete-time dynamical system\cite{Bertsekas1996} that follows the conventions of trajectory optimization.\cite{Lantoine2012p1, Lantoine2012p2}

%
%
\subsection{Discrete-Time Dynamical System and Control Policy}

Given a state vector $\bm{x}_k\in\mathbb{R}^{n_x}$ and a random vector $\bm{w}_k: \Omega\mapsto\mathbb{R}^{n_w}$, whose probability space is given as $(\Omega, \mathcal{F}, \mathbb{P})$, a discrete-time dynamical system with uncertainties is formulated by the stochastic equation
\begin{equation}
    \bm{x}_{k+1} = \bm{f}_k(\bm{x}_k, \bm{u}_k, \bm{w}_k), \ \ \ k \in \mathbb{N}_{N-1},\label{eq:dynamical_system}
\end{equation}
where $\bm{x}_{k+1}\in\mathbb{R}^{n_x}$ represents a successive state vector. The uncertainty vectors $\{\bm{w}_k\}_{k\in\mathbb{N}_{N-1}}$ are independent from each other, and the control vector $\bm{u}_k$ is the realization of a control policy $\bm{\mu}_k(\cdot)$.

%
%

This paper adopts closed-loop control to correct the perturbations that exist in stochastic systems. At each discrete time $k$, a control vector $\bm{u}_k$ is selected as the realization of a control policy $\bm{\mu}_k(\bm{x}_k)$ with the knowledge of the state vector $\bm{x}_k$. 
\begin{equation}
\bm{u}_k = \bm{\mu}_k(\bm{x}_k),  \ \ \ k \in \mathbb{N}_{N-1}.\label{eq:control_policy}
\end{equation}
Equation (\ref{eq:control_policy}) describes \textit{Markov control policy}\cite{Bertsekas1996}, which depends only on the current state.

%
%

\begin{figure}[htbp]
    \begin{center}
        \includegraphics[clip,width=0.65\hsize]{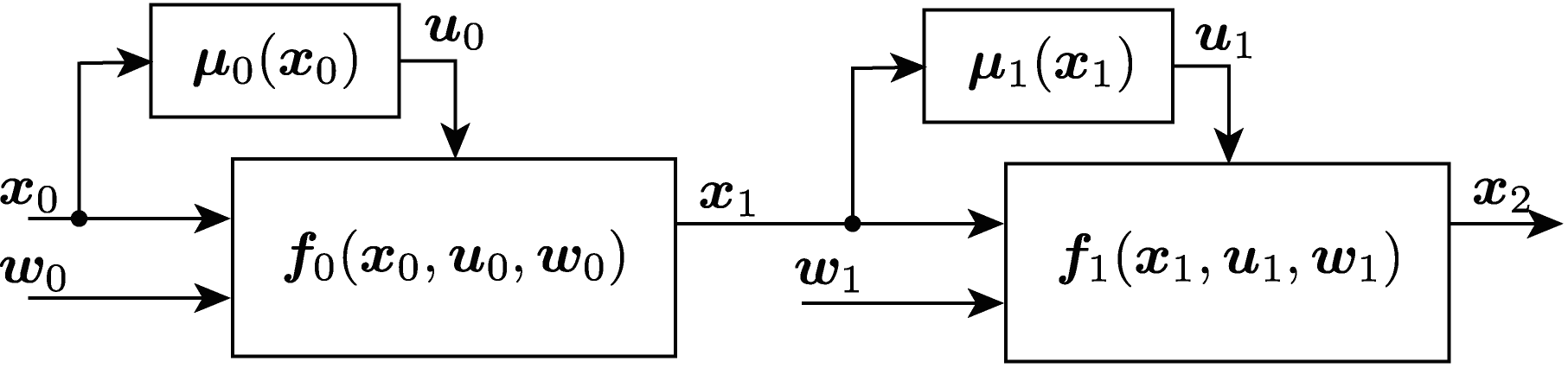}
        \caption{Propagation of dynamical system from $k=0$ to $k=2$.}
        \label{fig:propagation_of_dynsys}
    \end{center}
\end{figure}

This section describes the propagation process of the system illustrated as a block diagram in Fig. \ref{fig:propagation_of_dynsys}. The trajectory propagation starts with the selection of the state vector $\bm{x}_0$ from its probability space. The control vector $\bm{u}_0$ is decided through the control policy $\bm{\mu}_0(\bm{x}_0)$ with the state vector $\bm{x}_0$. The uncertainty vector $\bm{w}_0$ is also selected from its probability space. Finally, the state vector $\bm{x}_1$ is obtained using Eq.(\ref{eq:dynamical_system}). Continuing this propagation process up to $k=N$ yields the sequence of state vectors $\{\bm{x}_k\}_{k\in\mathbb{N}_{N}}$. Selecting a sample of $\bm{x}_0$ and $\{\bm{w}_k\}_{k\in\mathbb{N}_{N-1}}$ uniquely determines a trajectory $\{\bm{x}_k\}_{k\in\mathbb{N}_{N}}$. In actual stochastic systems, because the random variables $\bm{x}_0$ and $\{\bm{w}_k\}_{k\in\mathbb{N}_{N-1}}$ are not specified, the trajectory $\{\bm{x}_k\}_{k\in\mathbb{N}_{N}}$ is given as a stochastic process, that is, a set of random variables.

%
%
\subsection{Stochastic Optimal Control Problem}\label{sec:stochastic_optimal_control}
In a discrete-time system, the stochastic optimal control problem finds the set of the control policy $\bm{\pi} = \left\{\bm{\mu}_k(\cdot)\right\}_{k\in\mathbb{N}_{N-1}}$ that minimizes the objective function
\begin{equation}
J(\bm{x}_0,\bm{\pi}) = \mathbb{E} \left[\sum_{k=0}^{N-1} l_k(\bm{x}_k, \bm{u}_k, \bm{w}_k) + \varphi(\bm{x}_{N}) \right] \label{eq:soc_objective_function}
\end{equation}
subject to
\begin{align}
\bm{x}_{k+1} &= \bm{f}_k(\bm{x}_k, \bm{u}_k, \bm{w}_k), \ \ \ k\in\mathbb{N}_{N-1},\label{eq:soc_dynamical_system}\\
\bm{u}_k &= \bm{\mu}_k(\bm{x}_k),  \ \ \ k\in\mathbb{N}_{N-1},\label{eq:soc_control_policy}\\
\bm{\mu}_k(\bm{x}) &\in \mathbb{U}_k,  \forall \bm{x}\in \mathbb{X}_k \ \ \ k\in\mathbb{N}_{N-1}\label{eq:soc_control_constraints}
\end{align}
where the expected value $\mathbb{E}[\cdot]$ is taken with respect to the product of $\bm{x}_0$ and $\{\bm{w}_k\}_{k\in\mathbb{N}_{N-1}}$; $\mathbb{U}_k$ is the set of the admissible control; and $\mathbb{X}_k$ is the set of all possible $\bm{x}_k$.


%
%

Some specific problems, including linear quadratic ones\cite{Bertsekas2000}, can be solved analytically. On the other hand, generic nonlinear problems must be optimized numerically. However, implementing the numerical algorithm presents a number of challenges. 
The computation of the objective function (\ref{eq:soc_objective_function}) suffers from combinatorial explosion because the expected value must be taken with respect to the product of all random variables $\bm{x}_0$ and $\{\bm{w}_k\}_{k\in\mathbb{N}_{N-1}}$.
The propagation of the nonlinear dynamical system (\ref{eq:soc_dynamical_system}) with random variables causes difficulties because the high-order probabilistic moments are not negligible under the presence of nonlinearities.
The numerical expression of the control policy (\ref{eq:soc_control_policy}) requires further discussion because the numerical algorithm cannot search for the unknown mapping that describes the control policy, and it needs the approximation and parametrization of the control policy.
The handling algorithm of the control constraint (\ref{eq:soc_control_constraints}) cannot utilize conventional methods because the constraint on the possible control input associated with a random variable $\bm{x}$ is modeled by probabilistic constraints, which are more complex than deterministic ones. 
Finally, the implementation of a stable numerical optimization algorithm is also challenging because stochastic optimal control problems have more variables than conventional deterministic problems and involve more local optima, which should be avoided.

%% file: contents/3_tsoc_by_ut.tex
%
%
\section{Tube Stochastic Optimal Control by Unscented Transform}

This section presents a novel numerical method to solve the stochastic optimal control problem, as described in Eqs. (\ref{eq:soc_objective_function}) to (\ref{eq:soc_control_constraints}). This method approximates the state vector $\bm{x}_k$ with a Gaussian random variable at each discrete time $k$, employs the unscented transform to estimate the probability distribution of the state vectors, and parametrizes the control policies $\bm{\mu}_k(\cdot)$ by interpolating the control vectors on the sigma points of $\bm{x}_k$. These assumptions transform the stochastic optimal control problem into a deterministic optimal control problem.

%
%
\subsection{Assumptions}

This paper assumes that the initial state $\bm{x}_0$ and uncertainties $\bm{w}_k$ belong to Gaussian distributions $\mathcal{N}(\bar{\bm{x}}_0, \bm{\mathrm{P}}_0)$ and $\mathcal{N}(\bm{0}_{n_w}, \bm{\mathrm{R}}_k)$, respectively. At initial time $k=0$, the random state vector $\bm{x}_1$ can be calculated through the mappings $\bm{f}_0(\cdot)$ and $\bm{\mu}_0(\cdot)$ from $\bm{x}_{0}\sim\mathcal{N}(\bar{\bm{x}}_0,\bm{\mathrm{P}}_0)$ and $\bm{w}_{0}\sim\mathcal{N}(\bm{0}_{n_w},\bm{\mathrm{R}}_0)$ if the control policy $\bm{\mu}_0(\cdot)$ is given. The probability distribution of $\bm{x}_1$ becomes non-Gaussian because of the nonlinearity of $\bm{f}_0(\cdot)$ and $\bm{\mu}_0(\cdot)$, and its exact quantitative evaluation is impractical. Hence, this paper assumes that $\bm{x}_1$ belongs to a Gaussian distribution $\mathcal{N}(\bar{\bm{x}}_1, \bm{\mathrm{P}}_1)$, and evaluates $\bar{\bm{x}}_1$ and $\bm{\mathrm{P}}_1$ by the unscented transform. At any discrete time $k$, the random state vector $\bm{x}_{k+1}$ is calculated from $\bm{x}_{k}\sim\mathcal{N}(\bar{\bm{x}}_k,\bm{\mathrm{P}}_k)$ and $\bm{w}_{k}\sim\mathcal{N}(\bm{0}_{n_w},\bm{\mathrm{R}}_k)$, and its probability distribution is approximated as $\mathcal{N}(\bar{\bm{x}}_{k+1}, \bm{\mathrm{P}}_{k+1})$ by the unscented transform. The details of the unscented transform are summarized in Appendix A. Out of various uncertainty quantification methods, we select the unscented transform for its low computational burden and non-intrusive approach. Future work will adopt higher-order methods such as polynomial chaos expansion\cite{Xiu2003, Xiu2010, Jones2013}.

%
%
\subsection{Introduction of Sigma Points and Parametrization of Control Policies}

As the first step to transform the stochastic problem into a deterministic one, we introduce the sigma points of the state vector $\bm{x}_k$ and uncertainty vector $\bm{w}_k$. At each discrete time $k$, the probability distributions of $\bm{x}_k$ and $\bm{w}_k$ are given as the Gaussian distributions $\mathcal{N}(\bar{\bm{x}}_k, \bm{\mathrm{P}}_k)$ and $\mathcal{N}(\bm{0}_{n_w}, \bm{\mathrm{R}}_k)$, respectively. 

For the state vector $\bm{x}_k$, the unscented transform introduces the sigma points $\left\{\bm{\mathcal{X}}_k^{(i)}\right\}_{i\in\mathbb{N}_{2n_x}}$ and corresponding weights $\left\{c_x^{(i)}\right\}_{i\in\mathbb{N}_{2n_x}}$ as follows.
\begin{equation}
    \begin{cases}
        \bm{\mathcal{X}}_k^{(0)} &= \bar{\bm{x}}_k\\
        \bm{\mathcal{X}}_k^{(j)} &= \bar{\bm{x}}_k +  \left(\sqrt{(n_x+\kappa_x)\bm{\mathrm{P}}_k}\right)_j, \ \ j\in\mathbb{N}_{1:n_x}\\
        \bm{\mathcal{X}}_k^{(j+n_x)} &= \bar{\bm{x}}_k -  \left(\sqrt{(n_x+\kappa_x)\bm{\mathrm{P}}_k}\right)_j, \ \ j\in\mathbb{N}_{1:n_x}
    \end{cases},
\end{equation}
and
\begin{equation}
    c_x^{(i)} = \begin{cases}
        \frac{\kappa_x}{n_x + \kappa_x} &\text{     if } i=0\\
        \frac{1}{2(n_x+\kappa_x)} &\text{     if } i\in\mathbb{N}_{1:2n_x}
    \end{cases},
\end{equation}
where $\sqrt{\cdot}$ is the square root of the matrix and $(\cdot)_j$ is the $j$-th column vector of the matrix; and $\kappa_x$ is an arbitrary parameter for the unscented transform. Julier and Uhlmann\cite{Julier1996} have investigated the optimal choice of $\kappa_x$ to minimize the approximation errors.

Likewise, regarding the uncertainty vector $\bm{w}_k$, which is independent of the state vector $\bm{x}_k$, the sigma points $\left\{\bm{\mathcal{W}}_k^{(i)}\right\}_{i\in\mathbb{N}_{2n_w}}$ and corresponding weights $\left\{c_w^{(i)}\right\}_{i\in\mathbb{N}_{2n_w}}$ are defined as follows.

\begin{equation}
\begin{cases}
\bm{\mathcal{W}}_k^{(0)} &= \bm{0}_{n_w}\\
\bm{\mathcal{W}}_k^{(j)} &= \left(\sqrt{(n_w+\kappa_w)\bm{\mathrm{R}}_k}\right)_j, \ \ j\in\mathbb{N}_{1:n_w}\\
\bm{\mathcal{W}}_k^{(j+n_w)} &= -\left(\sqrt{(n_w+\kappa_w)\bm{\mathrm{R}}_k}\right)_j, \ \ j\in\mathbb{N}_{1:n_w}
\end{cases}
\end{equation}

and

\begin{equation}
c_w^{(i)} = \begin{cases}
\frac{\kappa_w}{n_w + \kappa_w} &\text{     if } i=0\\
\frac{1}{2(n_w+\kappa_w)} &\text{     if } i\in\mathbb{N}_{1:2n_w}
\end{cases}
\end{equation}
where $\kappa_w$ is an arbitrary parameter for the unscented transform.

In this paper, the sigma points of $\bm{x}_k$ and $\bm{w}_k$ are set up independently, and therefore the total number of sigma points is $(2n_x+1)\cdot(2n_w+1)$. An alternative approach defines a new random variable $\bm{y}_k:=[\bm{x}_k^T, \bm{w}_k^T]^T$ and creates the sigma points of $\bm{y}_k$. This approach reduces the total number of sigma points to $(2n_x+2n_w+1)$ and computational burden. However, independent sigma points, can capture the coupling term between $\bm{x}_k$ (or $\bm{u}_k$) and $\bm{w}_k$ more accurately than the alternative approach. Hence, we independently set up the sigma points of $\bm{x}_k$ and $\bm{w}_k$.

We now introduce the parameters representing the control policy (\ref{eq:soc_control_policy}). Given a control policy $\bm{\mu}_k(\cdot)$, the control vectors at the sigma points are expressed as
\begin{equation}
\bm{\mathcal{U}}_k^{(i)} = \bm{\mu}_k\left(\bm{\mathcal{X}}_k^{(i)}\right), \ \ \ i \in\mathbb{N}_{2n_x}.
\end{equation}
Although the control policy $\bm{\mu}_k(\cdot)$ should be defined for all possible state vectors, the unscented transform formulation only uses the control vectors at the sigma points $\left\{\bm{\mathcal{U}}_k^{(i)}\right\}_{i\in\mathbb{N}_{2n_x}}$. Thus, the proposed method finds the optimal control vectors $\left\{\bm{\mathcal{U}}_k^{(i)}\right\}_{i\in\mathbb{N}_{2n_x}}$ as parameters of the optimal control problem, and the actual control policy $\bm{\mu}_k(\cdot)$ is determined to be fitted to the sigma points $\left\{\bm{\mathcal{U}}_k^{(i)}\right\}_{i\in\mathbb{N}_{2n_x}}$.

Let $\bm{X}_k, \bm{U}_k,$ and $\bm{W}_k$ be the set of sigma points defined as follows:
\begin{align}
\bm{X}_k &:= \begin{bmatrix}
\bm{\mathcal{X}}^{(0)}_k & \bm{\mathcal{X}}^{(1)}_k & \cdots & \bm{\mathcal{X}}^{(2n_x)}_k
\end{bmatrix}\in\mathbb{R}^{n_x (2n_x+1)} \label{eq:set_of_sigma_points_Xk}\\
\bm{U}_k &:= \begin{bmatrix}
\bm{\mathcal{U}}^{(0)}_k & \bm{\mathcal{U}}^{(1)}_k & \cdots & \bm{\mathcal{U}}^{(2n_x)}_k
\end{bmatrix}\in\mathbb{R}^{n_u (2n_x+1)} \label{eq:set_of_sigma_points_Uk}\\
\bm{W}_k &:= \begin{bmatrix}
\bm{\mathcal{W}}^{(0)}_k & \bm{\mathcal{W}}^{(1)}_k & \cdots & \bm{\mathcal{W}}^{(2n_w)}_k
\end{bmatrix}\in\mathbb{R}^{n_w (2n_w+1)}.\label{eq:set_of_sigma_points_Wk}
\end{align}
The following sections derive the deterministic expression of the objective function (\ref{eq:soc_objective_function}), dynamical system (\ref{eq:soc_dynamical_system}), and control constraints (\ref{eq:soc_control_constraints}) with respect to $\bm{X}_k, \bm{U}_k,$ and $\bm{W}_k$. 

%
%
\subsection{Deterministic Expression of the Objective Function}

The objective function (\ref{eq:soc_objective_function}) can be expanded as
\begin{equation}
J(\bm{x}_0,\bm{\pi}) = \sum_{k=0}^{N-1} \underset{\bm{x}_k, \bm{w}_k}{\mathbb{E}}\left[l_k(\bm{x}_k, \bm{u}_k, \bm{w}_k)\right] + \underset{\bm{x}_N}{\mathbb{E}} \left[ \varphi(\bm{x}_{N})\right]\label{eq:soc_objective_function_2}
\end{equation}
because of the linearity of the expected value operator $\mathbb{E}$, i.e., $\mathbb{E}[x+y] = \mathbb{E}[x]+\mathbb{E}[y]$.

Using the set of sigma points $\bm{X}_k, \bm{U}_k,$ and $\bm{W}_k$ and taking a weighted sum of the representative value at the sigma points yields the expected values $\mathbb{E}\left[l_k(\bm{x}_k, \bm{u}_k, \bm{w}_k)\right]$ and $\mathbb{E} \left[ \varphi(\bm{x}_{N})\right]$ as follows:
\begin{align}
\underset{\bm{x}_k, \bm{w}_k}{\mathbb{E}}\left[l_k(\bm{x}_k, \bm{u}_k, \bm{w}_k)\right] &\simeq L_k(\bm{X}_k, \bm{U}_k, \bm{W}_k) := \sum_{i=0}^{2n_x} \sum_{j=0}^{2n_w} c_x^{(i)} c_w^{(j)} l_k\left(\bm{\mathcal{X}}_k^{(i)}, \bm{\mathcal{U}}_k^{(i)}, \bm{\mathcal{W}}_k^{(j)}\right), \label{eq:soc_objective_function_l}\\
\underset{\bm{x}_N}{\mathbb{E}} \left[ \varphi(\bm{x}_{N})\right] &\simeq \Phi(\bm{X}_{N}) := \sum_{i=0}^{2n_x} c_x^{(i)} \varphi \left(\bm{\mathcal{X}}_N^{(i)}\right)\label{eq:soc_objective_function_phi}
\end{align}
Substituting Eqs. (\ref{eq:soc_objective_function_l}) and (\ref{eq:soc_objective_function_phi}) into Eq. (\ref{eq:soc_objective_function_2}) gives the deterministic expression of the objective function
\begin{equation}
J_D(\bm{X}_0,\bm{\Pi}) := \sum_{k=0}^{N-1} L_k(\bm{X}_k, \bm{U}_k, \bm{W}_k) + \Phi(\bm{X}_{N}),\label{eq:tsoc_objective_function_0}
\end{equation}
where $\bm{\Pi}:=\{\bm{U}_k\}_{k\in\mathbb{N}_{N-1}}$ is the set of the optimization parameters.

%
%
\subsection{Deterministic Expression of the Dynamical System}

\begin{figure}[tbp]
    \begin{center}
        \includegraphics[clip,width=0.65\hsize]{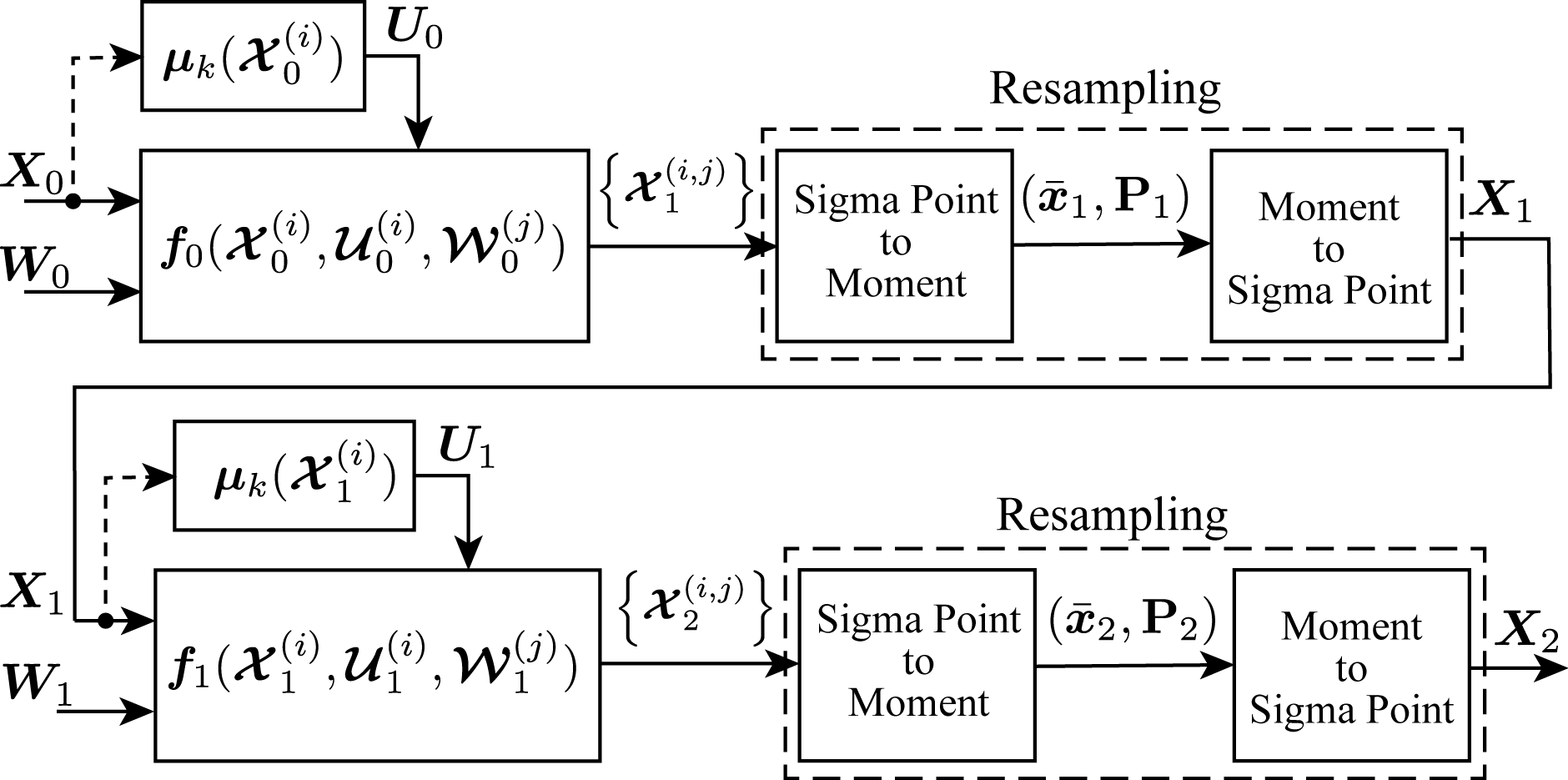}
        \caption{Block diagram of dynamical system propagation.}
        \label{fig:propagation_of_dynsys_ut}
    \end{center}
\end{figure}

This section derives the deterministic expression of the dynamical system, which has $(\bm{X}_k, \bm{U}_k, \bm{W}_k)$ as inputs and $\bm{X}_{k+1}$ as outputs, as shown in the block diagram in Fig. \ref{fig:propagation_of_dynsys_ut}. Given the sigma points $\bm{\mathcal{X}}_k^{(i)}, \bm{\mathcal{U}}_k^{(i)},$ and $\bm{\mathcal{W}}_k^{(j)}$ for $i\in\mathbb{N}_{2n_x}$ and $j \in\mathbb{N}_{2n_w}$ at discrete time $k$, the sigma points are propagated through the dynamical system (\ref{eq:soc_dynamical_system}) as follows.
\begin{equation}
\bm{\mathcal{X}}_{k+1}^{(i,j)} = \bm{f}_k \left(\bm{\mathcal{X}}_k^{(i)}, \bm{\mathcal{U}}_k^{(i)}, \bm{\mathcal{W}}_k^{(j)}\right), \ \ \ i\in\mathbb{N}_{2n_x}, j\in\mathbb{N}_{2n_w}\label{eq:dynamical_system_with_ut_step_1}
\end{equation}
Note that the dimension of the sigma points of the state vector increases exponentially at each step of the propagation, from $n_x\times(2n_x+1)$ to $n_x\times(2n_x+1)\times(2n_w+1)$. To limit the dimensionality of the problem, the proposed method resamples the sigma points at each step by approximating the probability distribution as a Gaussian distribution.

The first and second moments of the propagated state vector $\bm{x}_{k+1}$ are approximately computed as
\begin{align}
\bar{\bm{x}}_{k+1} &= \mathbb{E}\left[ \bm{f}_k(\bm{x}_k, \bm{u}_k, \bm{w}_k) \right]  \label{eq:dynamical_system_with_ut_step_2a0}\\
&\simeq \sum_{i=0}^{2n_x}\sum_{j=0}^{2n_w} c_x^{(i)} c_w^{(j)} \bm{\mathcal{X}}_{k+1}^{(i,j)},\label{eq:dynamical_system_with_ut_step_2a}\\
\bm{\mathrm{P}}_{k+1} &= \mathbb{E}\left[ \left\{ \bm{f}_k(\bm{x}_k, \bm{u}_k, \bm{w}_k) - \bar{\bm{x}}_{k} \right\} \left\{ \bm{f}_k(\bm{x}_k, \bm{u}_k, \bm{w}_k) - \bar{\bm{x}}_{k} \right\}^T \right]  \\
&\simeq \sum_{i=0}^{2n_x}\sum_{j=0}^{2n_w} c_x^{(i)} c_w^{(j)} \left\{\bm{\mathcal{X}}_{k+1}^{(i,j)} - \bar{\bm{x}}_{k+1}\right\} \left\{\bm{\mathcal{X}}_{k+1}^{(i,j)} - \bar{\bm{x}}_{k+1}\right\}^T\label{eq:dynamical_system_with_ut_step_2b}.
\end{align}
Finally, the resampled sigma points of $\bm{x}_{k+1}$ are computed as
\begin{align}
\bm{\mathcal{X}}_{k+1}^{(0)} &= \bar{\bm{x}}_{k+1},\label{eq:dynamical_system_with_ut_step_3a}\\
\bm{\mathcal{X}}_{k+1}^{(j)} &= \bar{\bm{x}}_{k+1} + \left(\sqrt{(n_x+\kappa)\bm{\mathrm{P}}_{k+1}}\right)_j,  \ \ \ j\in\mathbb{N}_{n_x}^+,\\
\bm{\mathcal{X}}_{k+1}^{(n+j)} &= \bar{\bm{x}}_{k+1} - \left(\sqrt{(n_x+\kappa)\bm{\mathrm{P}}_{k+1}}\right)_j, \ \ \ j\in\mathbb{N}_{n_x}^+.\label{eq:dynamical_system_with_ut_step_3c}
\end{align}
and the dimensions of the sigma points of the state vector remain $n_x\times(2n_x+1)$.

In summary, combining Eqs. (\ref{eq:dynamical_system_with_ut_step_1}) to (\ref{eq:dynamical_system_with_ut_step_3c}), we can formulate the propagation process by the deterministic nonlinear mapping $\bm{F}_k(\cdot)$
\begin{equation}
\bm{X}_{k+1} = \bm{F}_k(\bm{X}_k, \bm{U}_k, \bm{W}_k). \label{eq:dynamical_system_with_unscented_transform}
\end{equation}


%
%
\subsection{Deterministic Expression of Control Constraints}

\begin{figure}[tbp]
    \begin{center}
        \includegraphics[clip,width=0.65\hsize]{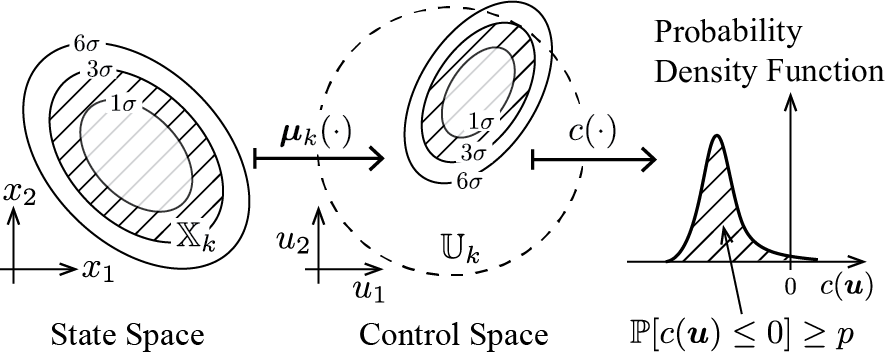}
        \caption{Chance-constrained method for control constraints.}
        \label{fig:control_constraints}
    \end{center}
\end{figure}

This section explains the procedure to formulate the deterministic expression of the control constraints (\ref{eq:soc_control_constraints}) under the assumption that they depend on control inputs only:
\begin{equation}
\mathbb{U}_k=\left\{\bm{u}\in\mathbb{R}^{n_u} : c(\bm{u})\leq 0\right\},
\end{equation}
implying
\begin{equation}
\forall \bm{x}_k \in \mathbb{X}_k, c\left(\bm{\mu}_k(\bm{x}_k)\right) \leq 0, \label{eq:control_deterministic_constraints}
\end{equation}
where $\mathbb{X}_k$ is the set of all the possible $\bm{x}_k$. $\mathbb{X}_k$ is equivalent to $\mathbb{R}^{n_x}$ whenever $\bm{x}_k$ is a Gaussian random variable. However, in this situation, the control constraint becomes numerically intractable due to the large number of possible realizations of $\bm{x}_k\in\mathbb{X}_k$.

To reduce the size of the domain, we introduce the chance-constrained method detailed in Appendix B, and rewrite Eq.(\ref{eq:control_deterministic_constraints}) as
\begin{equation}
\mathbb{P}\left[c\left(\bm{\mu}_k(\bm{x}_k)\right) \leq 0 : \bm{x}_k \sim \mathcal{N}(\bar{\bm{x}}_k, \bm{\mathrm{P}}_k)\right] \geq p, \label{eq:control_probabilistic_constraints}
\end{equation}
where $p\in[0,1]$ is the probability level. For example, $p=0.997$ for the 3-$\sigma$ satisfaction. Note that this constraint only guarantees the satisfaction of the original constraint (\ref{eq:control_deterministic_constraints}) in a probabilistic sense; however, this assumption is practically reasonable in trajectory optimization problems. Figure \ref{fig:control_constraints} offers an overview of the chance-constrained formulation.

Let us derive the deterministic expression of Eq.(\ref{eq:control_probabilistic_constraints}) using an example $\ell_2$-norm control constraint, which is among the most common constraint in low-thrust trajectory design. We can also extend the following calculation to general constraints $c(\bm{u})\leq 0$ in the same manner. Let us introduce the control constraint as
\begin{equation}
	c(\bm{\mu}_k(\bm{x}_k)) = \|\bm{\mu}_k(\bm{x}_k)\|^2 - u_{UB}^2 \leq 0
\end{equation}
and apply the unscented transform to the constraint. The first and second moments of $c(\bm{\mu}_k(\bm{x}_k))$ are calculated as
\begin{align}
\mathbb{E}\left[c(\bm{\mu}_k(\bm{x}_k))\right] &\simeq \sum_{i=0}^{2n_x} c_x^{(i)}\left\|\bm{\mu}_k(\bm{\mathcal{X}}_k^{(i)})\right\|^2 - u_{UB}^2 = \sum_{i=0}^{2n_x} c_x^{(i)}\left\|\bm{\mathcal{U}}_k^{(i)}\right\|^2 - u_{UB}^2,\\
\mathbb{V}\left[c(\bm{\mu}_k(\bm{x}_k))\right] &\simeq \sum_{i=0}^{2n_x} c_x^{(i)} \left\{\left\|\bm{\mu}_k(\bm{\mathcal{X}}_k^{(i)})\right\|^2 - u_{UB}^2 - \mathbb{E}\left[\left\|\bm{\mu}_k(\bm{x}_k)\right\|^2\right]\right\}^2\nonumber\\
&= \sum_{i=0}^{2n_x} c_x^{(i)} \left\{\left\|\bm{\mathcal{U}}_k^{(i)}\right\|^2 - \sum_{j=0}^{2n_x} c_x^{(j)}\left\|\bm{\mathcal{U}}_k^{(j)}\right\|^2 \right\}^2.
\end{align}
If we approximate $c(\bm{u})$ by a Gaussian distribution, Eq. (\ref{eq:control_probabilistic_constraints}) are written as
\begin{equation}
\mathbb{E}\left[c(\bm{\mu}_k(\bm{x}_k))\right] + 3\sqrt{\mathbb{V}\left[c(\bm{\mu}_k(\bm{x}_k))\right]} \leq 0,
\end{equation}

Therefore, the control constraint is described by a deterministic constraint using $\bm{U}_k$ as
\begin{equation}
C_k(\bm{U}_k) \leq 0,\label{eq:general_control_constraints}
\end{equation}
where
\begin{align}
C_k(\bm{U}_k) = \sum_{i=0}^{2n_x} c_x^{(i)}\left\|\bm{\mathcal{U}}_k^{(i)}\right\|^2 + 3\sqrt{\sum_{i=0}^{2n_x} c_x^{(i)} \left\{\left\|\bm{\mathcal{U}}_k^{(i)}\right\|^2 - \sum_{j=0}^{2n_x} c_x^{(j)}\left\|\bm{\mathcal{U}}_k^{(j)}\right\|^2 \right\}^2} - u_{UB}^{2}
\end{align}

Constrained optimal control problems with the constraints (\ref{eq:general_control_constraints}) have a singular point when the variance $\mathbb{V}[\cdot]$ approaches zero. To avoid this singularity during numerical optimization, this paper implements the following constraints
\begin{equation}
C_k(\bm{U}_k) = \mathbb{E}\left[c(\bm{\mu}_k(\bm{x}_k))\right] + 3\sqrt{\mathbb{V}\left[c(\bm{\mu}_k(\bm{x}_k))\right]+\epsilon} - 3\sqrt{\epsilon}\leq 0,\label{eq:tsoc_constraint_mod}
\end{equation}
where $\epsilon<<1$ is a small number to avoid the singularity ($\epsilon = 10^{-4}$ in our implementation).

%
%
\subsection{Tube Stochastic Optimal Control by the Unscented Transform}

The original stochastic optimal control problem (\ref{eq:soc_objective_function}) to (\ref{eq:soc_control_constraints}) is transformed into the deterministic problem, which finds the set of control vectors $\bm{\Pi}:=\left\{\bm{U}_k\right\}_{k\in\mathbb{N}_{N-1}}$ to minimize
\begin{equation}
J_D(\bm{X}_0,\bm{\Pi}) = \sum_{k=0}^{N-1} L_k(\bm{X}_k, \bm{U}_k, \bm{W}_k) + \Phi(\bm{X}_{N})\label{eq:tsoc_begin}
\end{equation}
subject to
\begin{align}
\bm{X}_{k+1} &= \bm{F}_k(\bm{X}_k, \bm{U}_k, \bm{W}_k), \ \ \ k\in\mathbb{N}_{N-1},\\
C_k(\bm{U}_k) &\leq 0, \ \ \ k\in\mathbb{N}_{N-1}\\
\bm{X}_0 &= \bar{\bm{X}}_0\label{eq:tsoc_end}
\end{align}
where $\bar{\bm{X}}_0$ and $\bm{W}_k$ are given from $\mathcal{N}(\bar{\bm{x}}_0, \bm{\mathrm{P}}_0)$ and $\mathcal{N}(\bm{0}_{n_w}, \bm{\mathrm{R}}_k)$; $J_D(\cdot)$ is given in Eq. (\ref{eq:tsoc_objective_function_0}); $\bm{F}_k(\cdot)$ is provided in Eq. (\ref{eq:dynamical_system_with_unscented_transform}); and $C_k(\cdot)$ is defined in Eq. (\ref{eq:tsoc_constraint_mod}). This deterministic optimal control problem can be solved by conventional trajectory optimization methods such as DDP\cite{Jacobson1970, Lantoine2012p1}.


%
%
\subsection{Differential Dynamic Programming}

DDP is among the most numerically stable methods that solves a second-order expansion of the Bellman equation to find the locally optimal solution. DDP has been widely used in space mission design\cite{Whiffen2002, Lantoine2012p1, Lantoine2012p2, Colombo2009, Pellegrini2012, Aziz2017} because it is suitable for large-scale optimal control problems and the solution converges to the optimal point robustly. This paper employs DDP to solve the optimal control problem (\ref{eq:tsoc_begin})-(\ref{eq:tsoc_end}) especially because of its numerical robustness, and calls the algorithm Tube Stochastic Differential Dynamic Programming (TSDDP).

%% file: contents/4_num_example_double_integ.tex
%
%
\section{Numerical Example 1: Double Integrator Problem}

In this section, a simple example is used to show that the solution computed by the proposed method is indeed both robust and optimal. The example problem is a 1-dimensional transfer between fixed initial and final states (positions and velocities) in fixed time, using bounded control inputs. The objective function is the sum of the integral of the control norm ($\Delta V$) and a penalty function that evaluates the violation of the final state constraint.

\subsection{Statement of the Problem}

For a state vector $\bm{x}_k=\left[r_k, v_k\right]^T\in\mathbb{R}^2$, a control vector $u_k\in[-1, 1]\subset\mathbb{R}$, and an uncertainty vector $\bm{w}_k\sim \mathcal{N}(\bm{0}, \bm{\mathrm{R}}_k)$, the equation of motion of a 1-dimensional double integrator is
\begin{equation}
\begin{bmatrix}
r_{k+1} \\ v_{k+1}
\end{bmatrix} = \begin{bmatrix}
1 & \Delta t\\
0 & 1
\end{bmatrix}\begin{bmatrix}
r_{k} \\ v_{k}
\end{bmatrix} + \begin{bmatrix}
0\\ b
\end{bmatrix} u_k + \bm{w}_k.\label{eq:dynsys_double_integ}
\end{equation}
where $\Delta t$ and $b$ are parameters that correspond to a time step and acceleration, respectively. The initial condition is given as $\bm{x}_0 \sim \mathcal{N}(\bar{\bm{x}}_0, O_2)$.

The objective function of the problem is formulated as the sum of the integral of the control norm and the violation of the final state constraint as follows:
\begin{equation}
J = \sum_{k=0}^{N-1} \mathbb{E}[\|u_k\|] + c_f \mathbb{E}\left[(\bm{x}_{N}-\bar{\bm{x}}_{N})^2\right]
\end{equation}
The parameters of the problem are summarized in Table \ref{tab:parameter_settings}.

\begin{table}[htb]
	\begin{center}
		\caption{Parameter settings}
		\begin{tabular}{lcc} \hline \hline
			Parameters & Variables & Settings \\ \hline
			Stage number & $N$ & 39 \\
			Time step & $\Delta t$ & 0.15\\
			Acceleration magnitude & $b$ & 0.25\\
			Covariance of uncertainties & $\bm{\mathrm{R}}_k$ & $\begin{bmatrix}
			10^{-20} & 0\\ 0 & 2.5\times 10^{-4}
			\end{bmatrix}$ \\ 
			Initial state & $\bar{\bm{x}}_0$ & $[-10, 0]^T$\\
			Final state & $\bar{\bm{x}}_N$ & $[0, 0]^T$\\ 
			Weight of objective function & $c_f$ & $10^4$\\ \hline \hline
		\end{tabular}
		\label{tab:parameter_settings}
	\end{center}
\end{table}

\subsection{Nominal Trajectories and Control}
\begin{figure}[htbp]
  \begin{center}
    \begin{tabular}{c}
	\hspace*{-0.075\hsize}
      \begin{minipage}{0.55\hsize}
		\begin{center}
			\includegraphics[clip,width=\hsize]{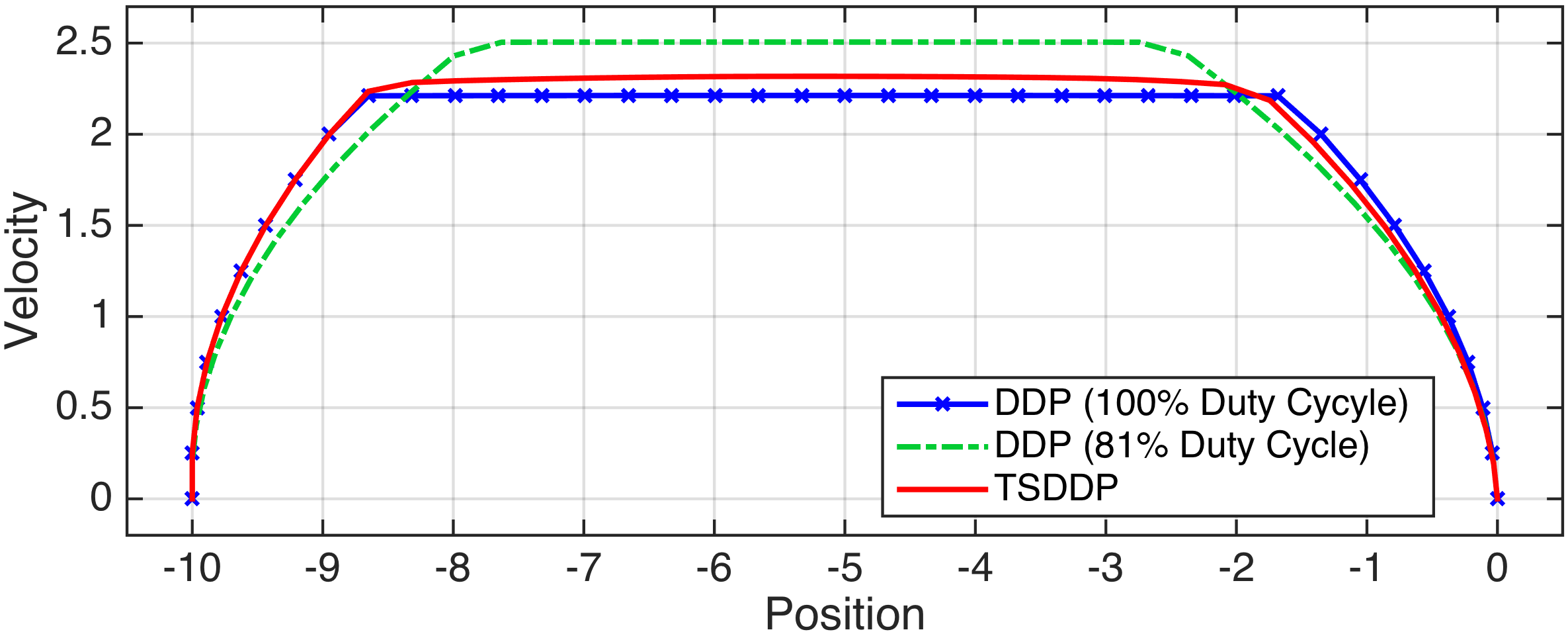}
			{\small a) Nominal trajectories in state space}
		\end{center}
      \end{minipage}

      \begin{minipage}{0.55\hsize}
		\begin{center}
			\includegraphics[clip,width=\hsize]{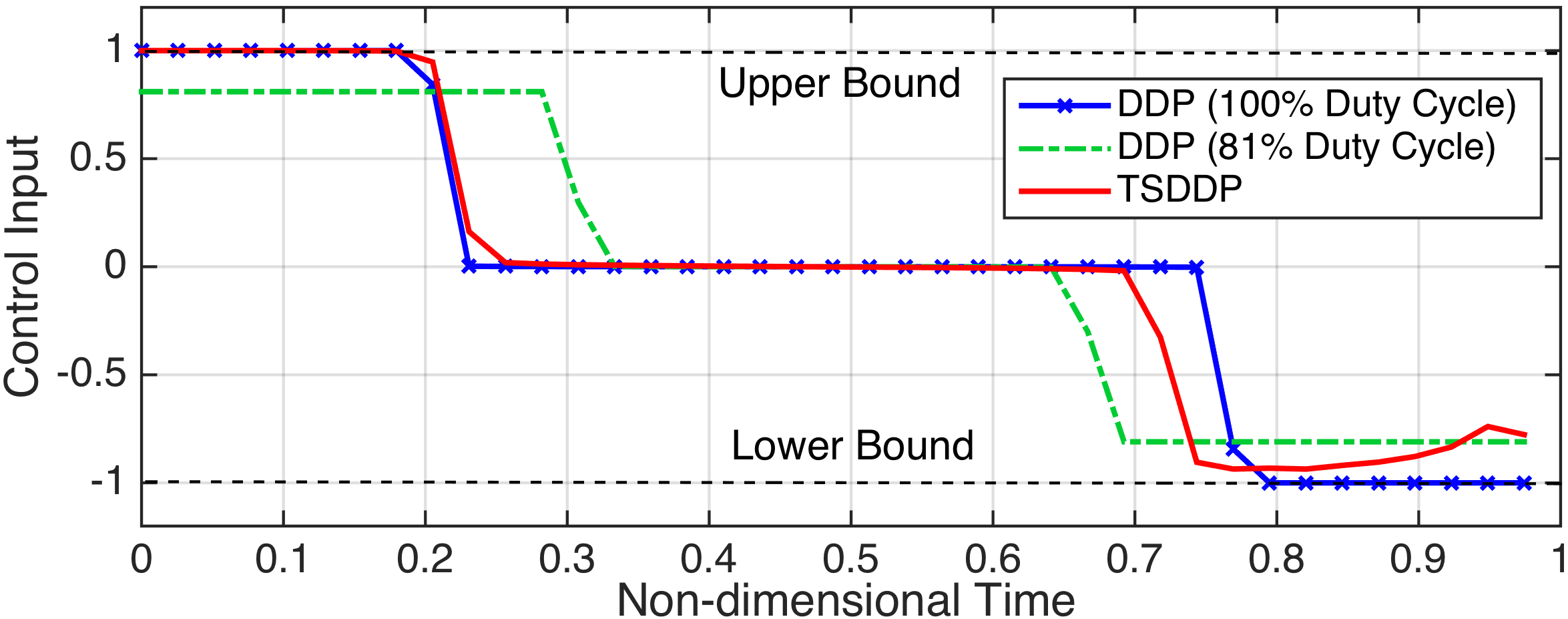}
			{\small b) Nominal control inputs}
		\end{center}
      \end{minipage}
    \end{tabular}
    \caption{Nominal trajectory and control inputs.}
    \label{fig:nom_trjcp_all}
  \end{center}
\end{figure}

Figure \ref{fig:nom_trjcp_all} shows the nominal trajectories and control computed by the conventional methods (DDP) and the proposed method (TSDDP). The conventional methods introduce duty cycles that reduce the the upper bound of the control input to improve robustness against uncertainties. The 100\% duty cycle case gives the less robust but most $\Delta V$-efficient solution, and the 81\% duty cycle case, where 81\% duty cycle is the required amount to compensate the $3\sigma$ uncertainty $\bm{w}_k$, gives a more robust but less $\Delta V$-efficient solution. The conventional solutions yield bang-bang control; contrarily, the proposed solution reduces the control input in the second thrusting arc to improve robustness, that is, the duty cycle. The results imply that the proposed method provides the appropriate amount of duty cycle automatically.

\subsection{Monte Carlo simulations}
\begin{figure}[htbp]
	\begin{center}
		\begin{tabular}{cc}
			\hspace*{-0.075\hsize}
			\begin{minipage}{0.55\hsize}
				\begin{center}
					\includegraphics[clip,width=\hsize]{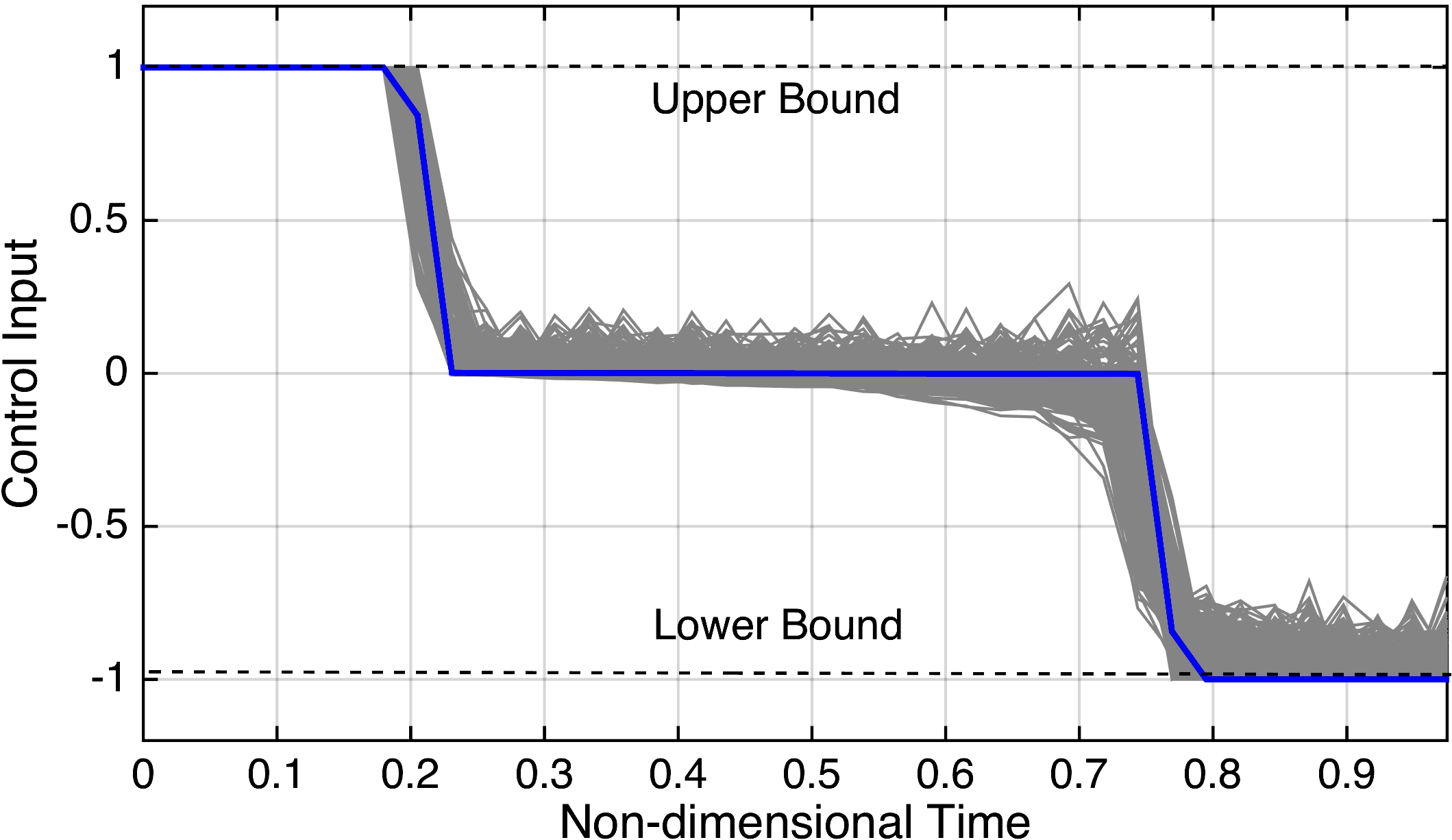}
					{\small a) DDP with 100\% duty cycle (case 1)}
					\vspace{7pt}
				\end{center}
			\end{minipage}
			\begin{minipage}{0.55\hsize}
				\begin{center}
					\includegraphics[clip,width=\hsize]{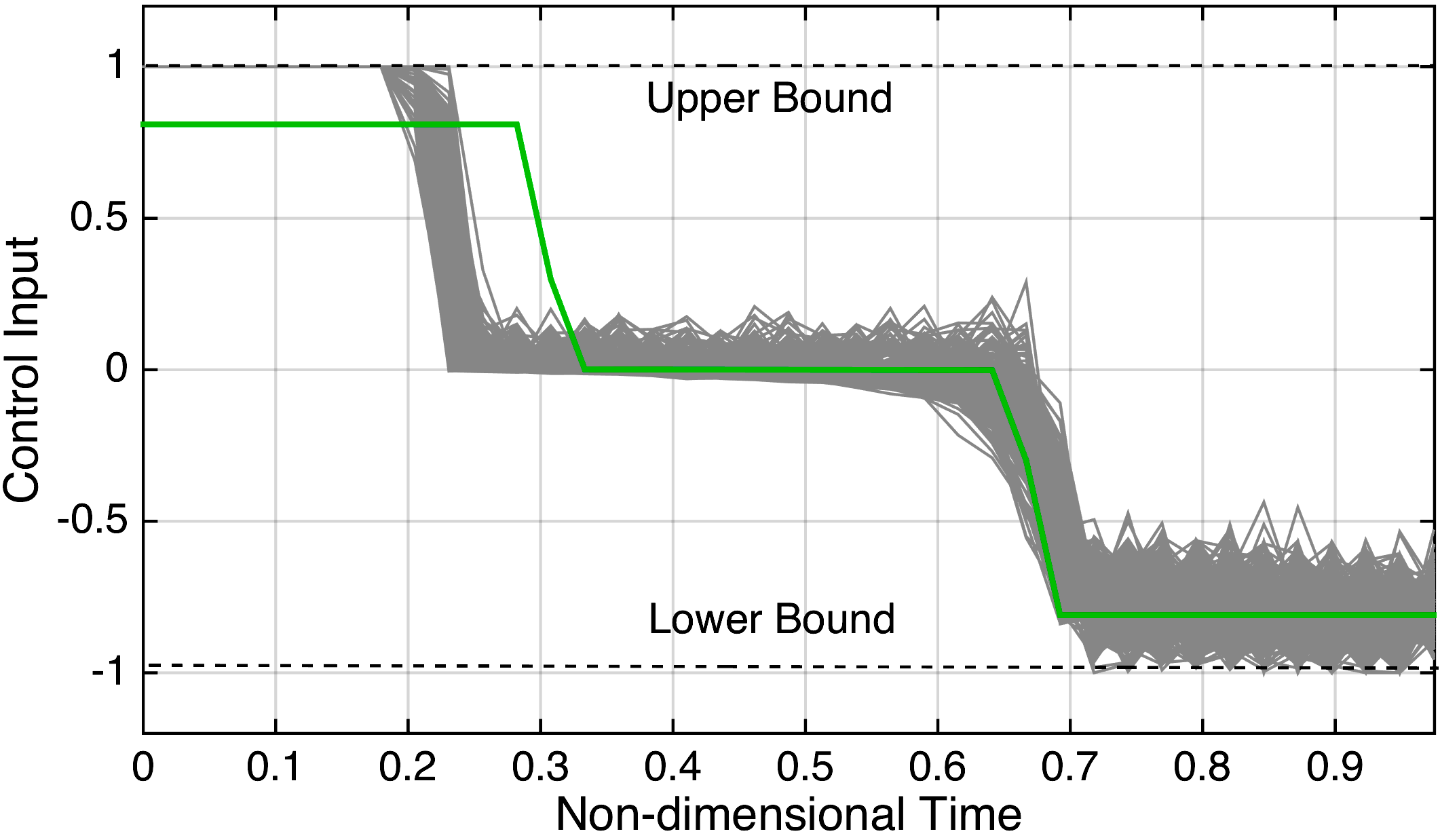}
					{\small b) DDP with 81\% duty cycle (case 2)}
					\vspace{7pt}
				\end{center}
			\end{minipage}
			\\
			\hspace*{-0.075\hsize}
			\begin{minipage}{0.55\hsize}
				\begin{center}
					\includegraphics[clip,width=\hsize]{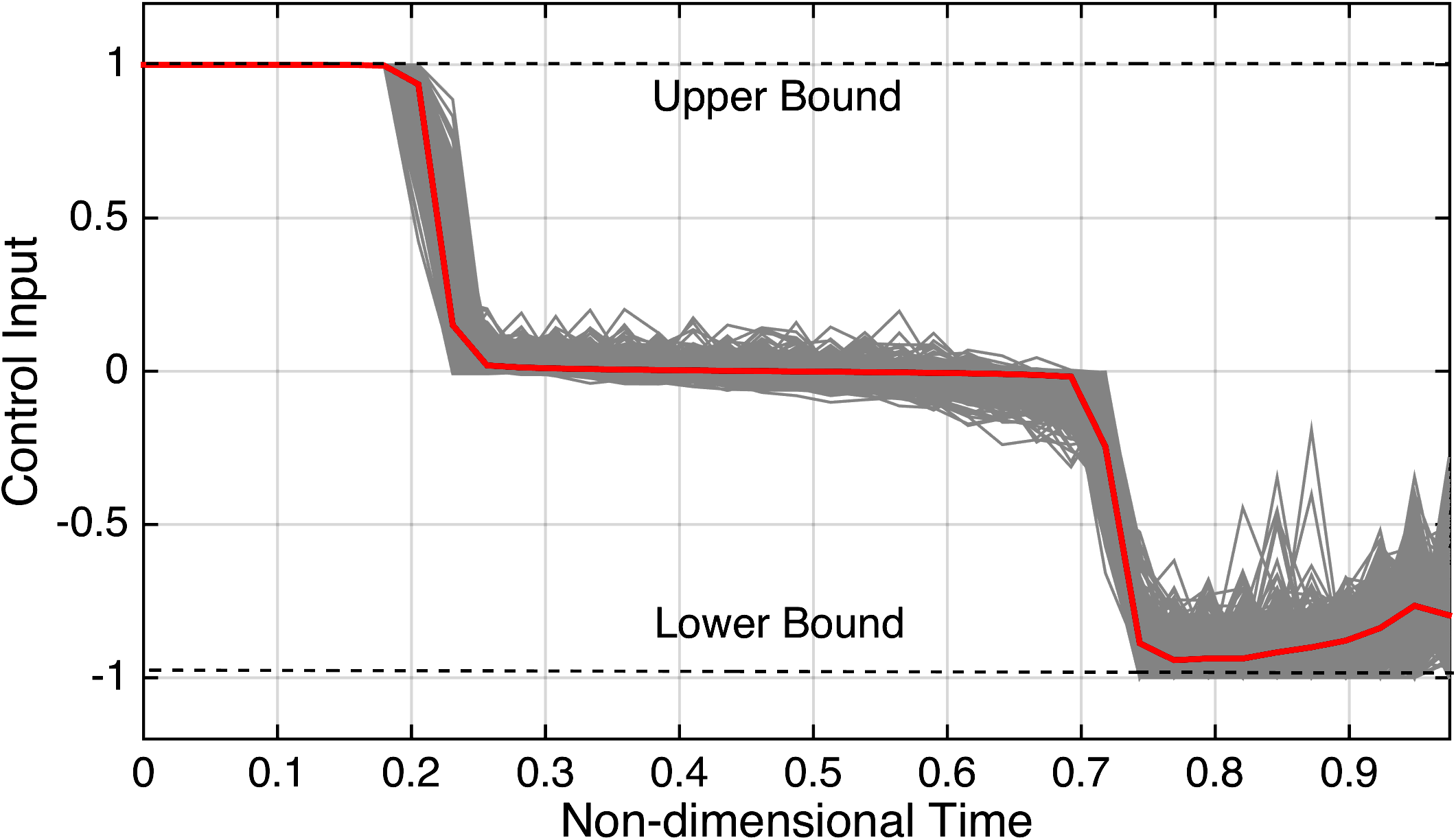}
					{\small c) TSDDP with re-optimization (case 3)}
				\end{center}
			\end{minipage}
			\begin{minipage}{0.55\hsize}
				\begin{center}
					\includegraphics[clip,width=\hsize]{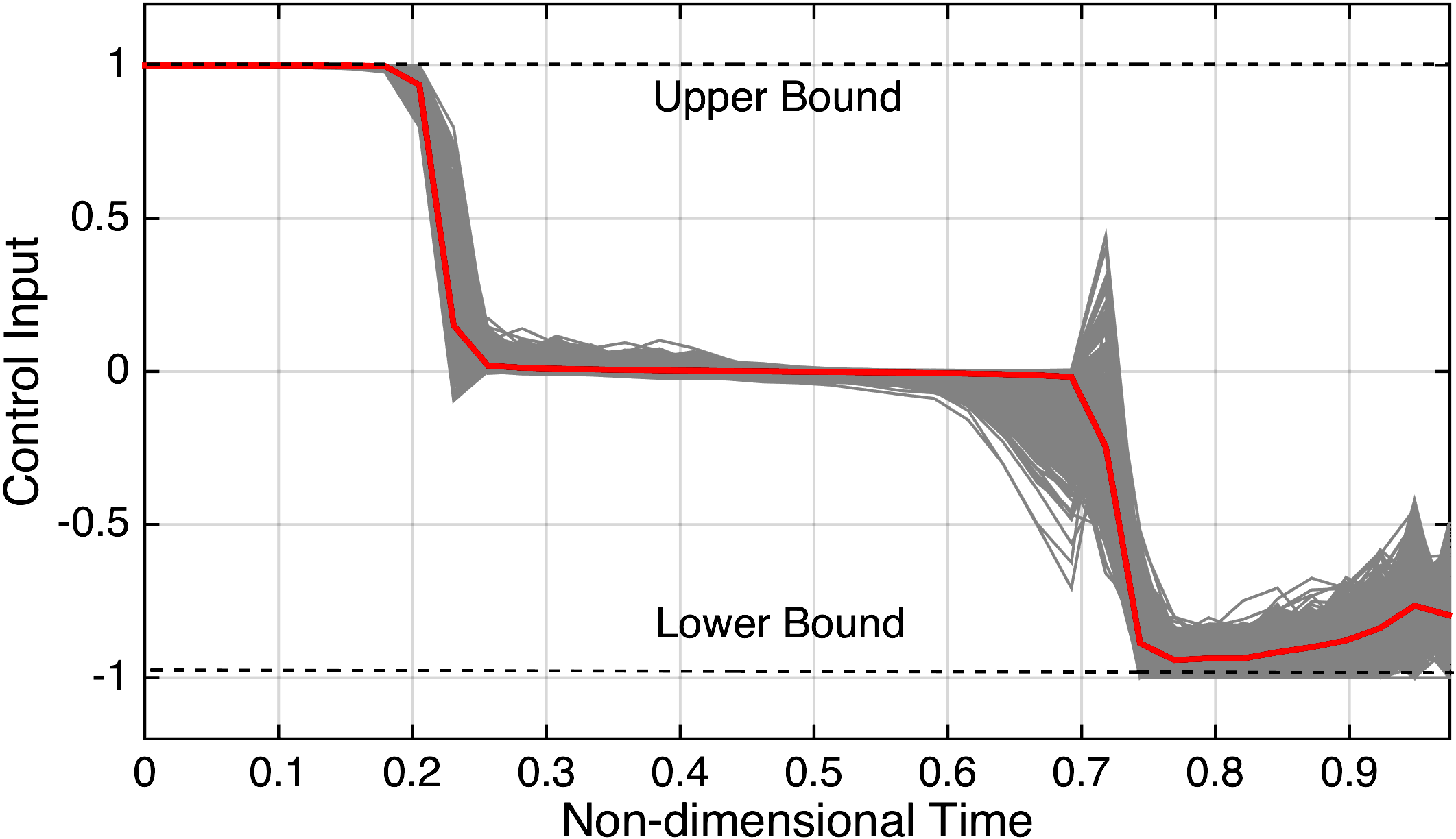}
					{\small d) TSDDP with interpolated control policy (case 3')}
				\end{center}
			\end{minipage}
		\end{tabular}
		\caption{Control profiles of Monte Carlo simulation. (500 samples, each gray line indicates a sample control profile and solid line represents the nominal control profile.)}
		\label{fig:mc_cp_all}
	\end{center}
\end{figure}

In this section, Monte Carlo simulations are run to observe the robustness and optimality of the proposed method. In each simulation, the optimization problems are recursively solved in a receding horizon fashion. Each sample trajectory of the Monte Carlo simulation is a numerically propagated solution of Eq.(\ref{eq:dynsys_double_integ}) from the initial discrete time ($k=0$) to the final discrete time ($k=N$). At each time step $k$, errors $\bm{w}_k$ are introduced and the remaining part of the trajectory is re-optimized using DDP with 100\% duty cycle (case 1), DDP with 81\% duty cycle (case 2), or TSDDP (case 3). Particularly, in the 81\% duty cycle case, each re-optimization at step $k$ sets 100\% duty cycle at the first segment (step $k$) because the state is fully given, and imposes 81\% duty cycle on the rest of the steps (step $(k+1)$ and later) as margins. Figure \ref{fig:mc_cp_all} plots the control profiles of the Monte Carlo simulation. For reference, Fig.\ref{fig:mc_cp_all}-d) shows the control profiles that are not re-optimized but use the piecewise linear control policy interpolating the control vectors at the sigma points.

\begin{figure}[htbp]
	\begin{center}
		\begin{tabular}{c}
			\begin{minipage}{0.6\hsize}
				\begin{center}
					\includegraphics[clip,width=\hsize]{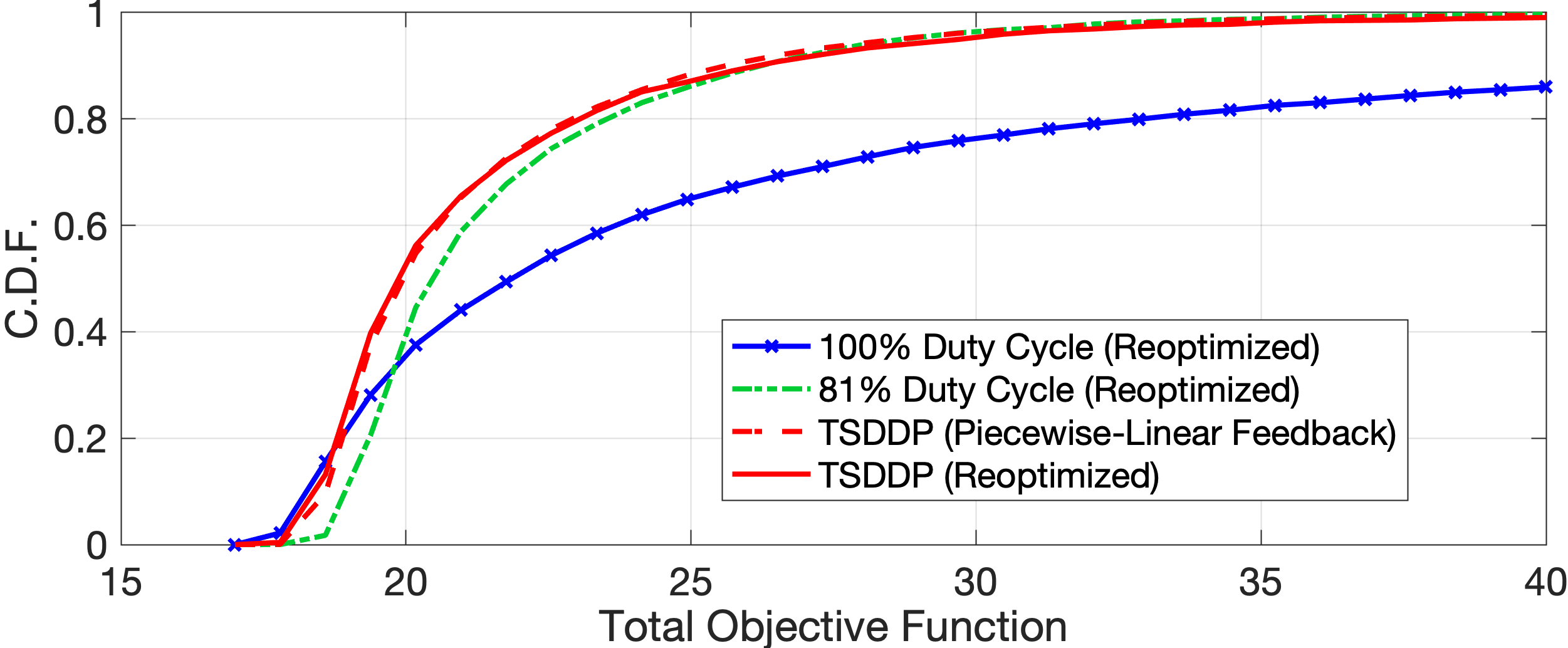}
				\end{center}
			\end{minipage}
			\\
			\begin{minipage}{0.6\hsize}
				\begin{center}
					\includegraphics[clip,width=\hsize]{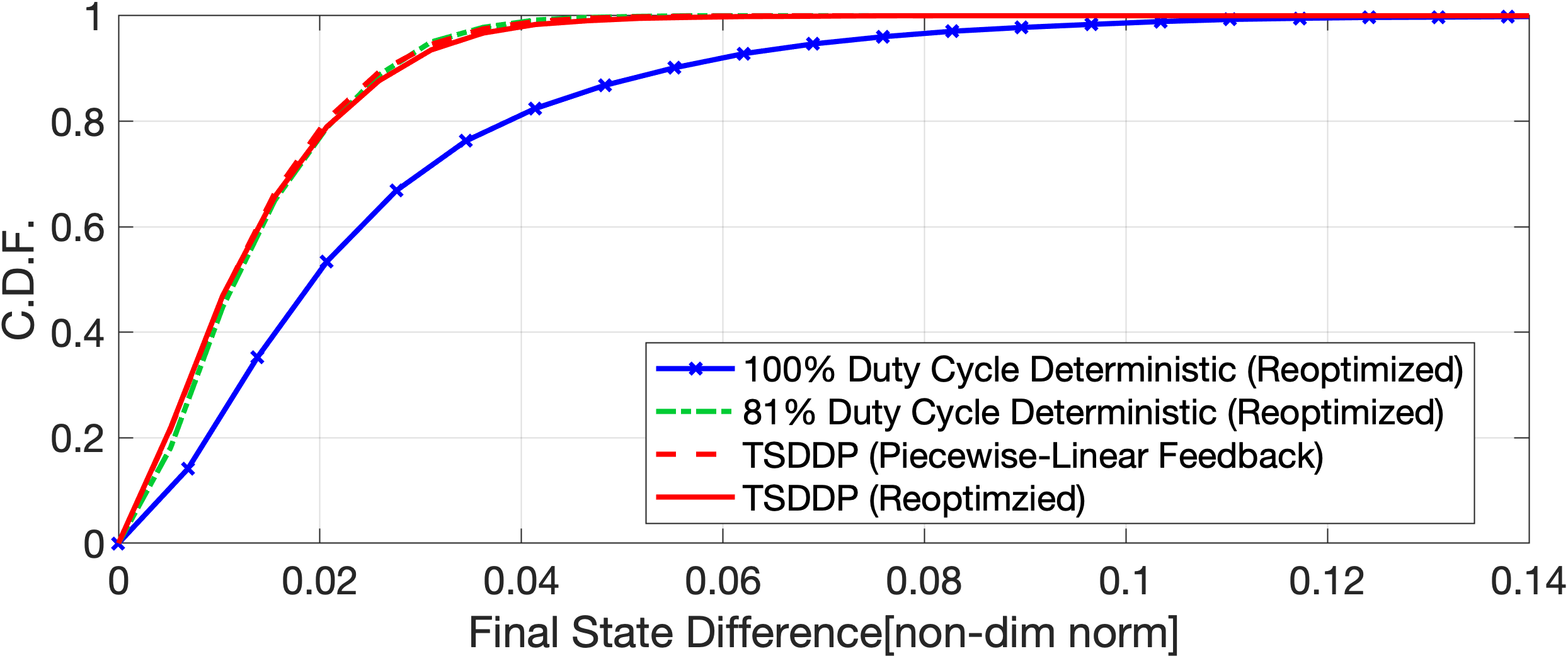}
				\end{center}
			\end{minipage}
			\\
			\begin{minipage}{0.6\hsize}
				\begin{center}
					\includegraphics[clip,width=\hsize]{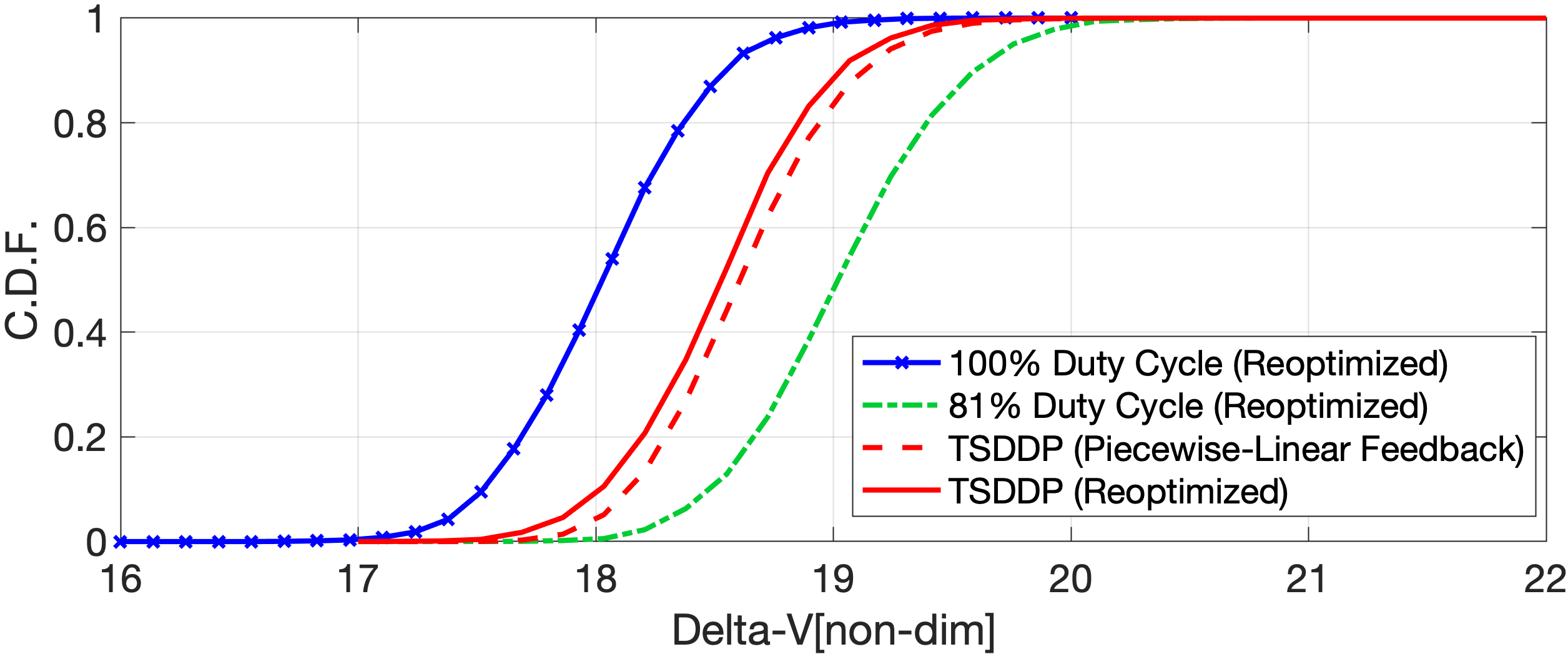}
				\end{center}
			\end{minipage}
		\end{tabular}
		\caption{Cumulative distribution functions (CDFs) of the objective functions.}
		\label{fig:mc_cdf_81p}
	\end{center}
\end{figure}

The robustness and efficiency of the solution are evaluated by the Cumulative Distribution Functions (CDFs) of the objective functions as plotted in Fig. \ref{fig:mc_cdf_81p}. 
The CDFs of the total objective functions show that TSDDP achieves the lowest cost in the three cases. Let us inspect the CDFs of two individual objectives, i.e., the final state difference and $\Delta V$, that balance each other. The CDFs of the final state difference, which indicate robustness against uncertainties, show that both case 2 and 3 are more robust than case 1. In addition, the CDFs of the $\Delta V$, which indicate efficiency, show that case 3 is more efficient than case 2.
Therefore, it is concluded that the proposed method is more robust and optimal than the conventional methods. Note that the conventional methods provide better results by a manual tuning of the duty cycle (parameters); nevertheless, the proposed method is useful since it gives robust and optimal solution without the need to tune these parameters.




%% file: contents/5_num_example_low_thrust.tex
%
%
\section{Numerical Example 2: Low-Thrust Trajectory Design}

In this section, the proposed method is applied to low-thrust trajectory design to show the robustness and efficiency in a nonlinear constrained system. The example deals with an Earth-Mars transfer problem in planar two-body dynamics in fixed time for a spacecraft equipped with low-thrust capabilities. In the same way as the first example, the objective function is the sum of the total $\Delta V$ magnitude and a penalty function that evaluates the violation of the final state constraint.

\subsection{Statement of the Problem}

For a state vector $\bm{x}=[\bm{r}^T, \bm{v}^T]^T\in \mathbb{R}^4$ and a control vector $\bm{u}\in\mathbb{U}\subset\mathbb{R}^2$, where $\mathbb{U}=\{\bm{u}\in\mathbb{R}^2: \|\bm{u}\|\leq u_{UB} \}$, the equations of motion of the spacecraft are
\begin{equation}
\frac{d}{dt}\begin{bmatrix}
\bm{r}\\ \bm{v}
\end{bmatrix} = \begin{bmatrix}
\bm{v}\\ -GM_{\odot} \bm{r}/\|\bm{r}\|^3
\end{bmatrix} + \begin{bmatrix}
\bm{0}_2\\ \bm{u}
\end{bmatrix}\label{eq:low_thrust_continuous_system}
\end{equation}
where $GM_{\odot}$ is the gravity constant of the Sun. 

Integrating Eq. (\ref{eq:low_thrust_continuous_system}) between $[t_k, t_{k+1})$ and adding an uncertainty vector $\bm{w}_k\sim \mathcal{N}(\bm{0}_4, \bm{\mathrm{R}}_k)$, where
\begin{equation}
\bm{\mathrm{R}}_k=\begin{bmatrix}\sigma_r^2 I_2 & O_2\\ O_2 & \sigma_v^2 I_2\end{bmatrix},
\end{equation}
we obtain the discrete-time stochastic dynamical system as
\begin{equation}
\bm{x}_{k+1} = \bm{f}_k (\bm{x}_k, \bm{u}_k) + \bm{w}_k,\label{eq:low_thrust_dynsys}
\end{equation}
This paper adopts the Runge-Kutta 4th-order method where the control vector $\bm{u}$ is kept constant between each discrete time step and $\bm{u}_k$ is the control vector between $[t_k, t_{k+1})$. The initial condition is given as
\begin{equation}
\bm{x}_0 \sim \mathcal{N}([\bm{r}_{\oplus}^T, \bm{v}_{\oplus}^T]^T, O_4).
\end{equation}

The objective function of the problem is formulated as the sum of the total $\Delta V$ and the violation of the final state constraint:
\begin{equation}
J = \sum_{k=0}^{N-1} \mathbb{E}\left[\sqrt{\|\bm{u}_k\|^2+\epsilon}\right] + c_f \mathbb{E}\left[(\bm{r}_{N}-\bm{r}_{\mars})^2 + (\bm{v}_{N}-\bm{v}_{\mars})^2\right],	
\end{equation}
where $\epsilon<<1$ is a mass-leak term used to prevent singularity in the computation of the partial derivatives. In our implementation, $\epsilon = 10^{-6}$. These parameters are summarized in Table \ref{tab:parameter_settings_low_thrust}. The values for the position and velocity of the Earth and Mars are the same as in Lantoine and Russell\cite{Lantoine2012p2}, but the orbits are projected onto the ecliptic plane.

\begin{table}[htb]
	\begin{center}
		\caption{Parameter settings}
		\begin{tabular}{lcc} \hline \hline
			Parameters & Variables & Settings \\ \hline
			Stage number & $N$ & 40 \\
			Transfer time (days) & $T$ & 348.79\\
			Gravity constant of the Sun & $GM_{\odot}$ & $1.32712442099\times10^{11}$\\
			Maximum thrust acceleration (km/s$^2$) & $u_{UB}$ & $1.0\times10^{-6}$\\ 
			Position errors (km$^2$) & $\sigma_r^2$ & $1.0\times 10^{-12}$\\
			Velocity errors (km$^2$/s$^2$) & $\sigma_v^2$ & $2.522627\times 10^{-5}$\\
			Initial position (km) & $\bm{r}_{\oplus}$ & $\begin{bmatrix} -140699693 & -51614428 \end{bmatrix}^T$ \\
			Initial velocity (km/s) & $\bm{v}_{\oplus}$ & $\begin{bmatrix} 9.774596 & -28.07828 \end{bmatrix}^T$ \\
			Final position (km) & $\bm{r}_{\mars}$ & $\begin{bmatrix} -172682023 & 176959469 \end{bmatrix}^T$ \\
			Final velocity (km/s) & $\bm{v}_{\mars}$ & $\begin{bmatrix} -16.427384 & -14.860506 \end{bmatrix}^T$ \\
			Weight of objective function & $c_f$ & $10^6$\\
			Scale factor & $L_{sf}, T_{sf}$ & $10^8, 10^6$\\
			Parameter of unscented transform & $\kappa$ & 2\\\hline \hline
		\end{tabular}
		\label{tab:parameter_settings_low_thrust}
	\end{center}
\end{table}

\subsection{Results and Discussion}
This section shows the nominal results and Monte Carlo analyses in the same manner as the previous example. We compare the proposed method (TSDDP) with the conventional methods (DDP with 100\% duty cycle and 80\% duty cycle), where 80\% is the maximum duty cycle to compensate the 3$\sigma$ value of the uncertainty $\bm{w}_k$. 

\begin{figure}[htbp]
  \begin{center}
    \begin{tabular}{c}
	\hspace*{-0.075\hsize}
      \begin{minipage}{0.55\hsize}
		\begin{center}
			\includegraphics[clip,width=\hsize]{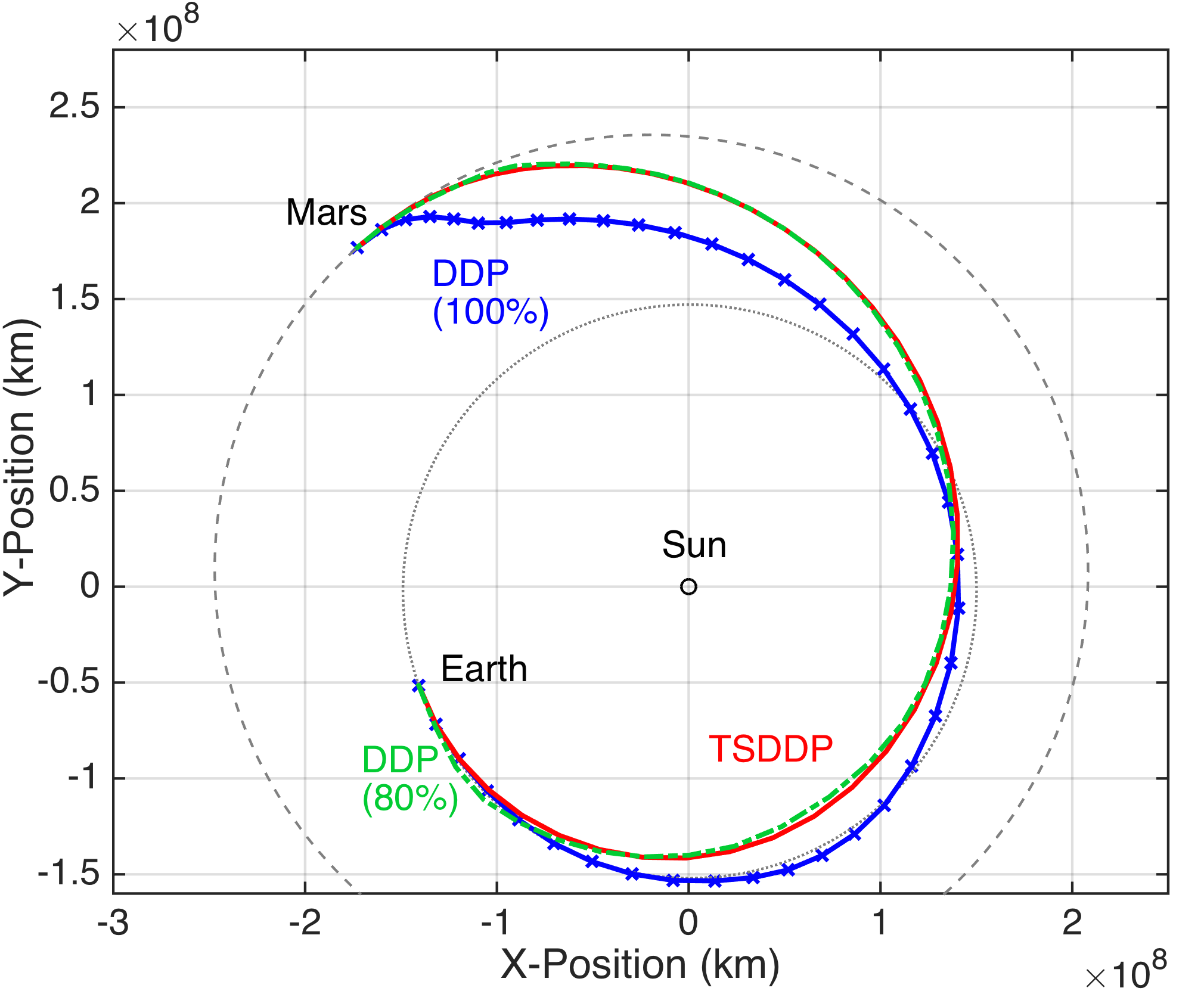}
			{\small a) Nominal trajectories (Differences from TSDDP are exaggerated 20 times for illustration purposes)}
		\end{center}
      \end{minipage}

      \begin{minipage}{0.55\hsize}
		\begin{center}
			\includegraphics[clip,width=\hsize]{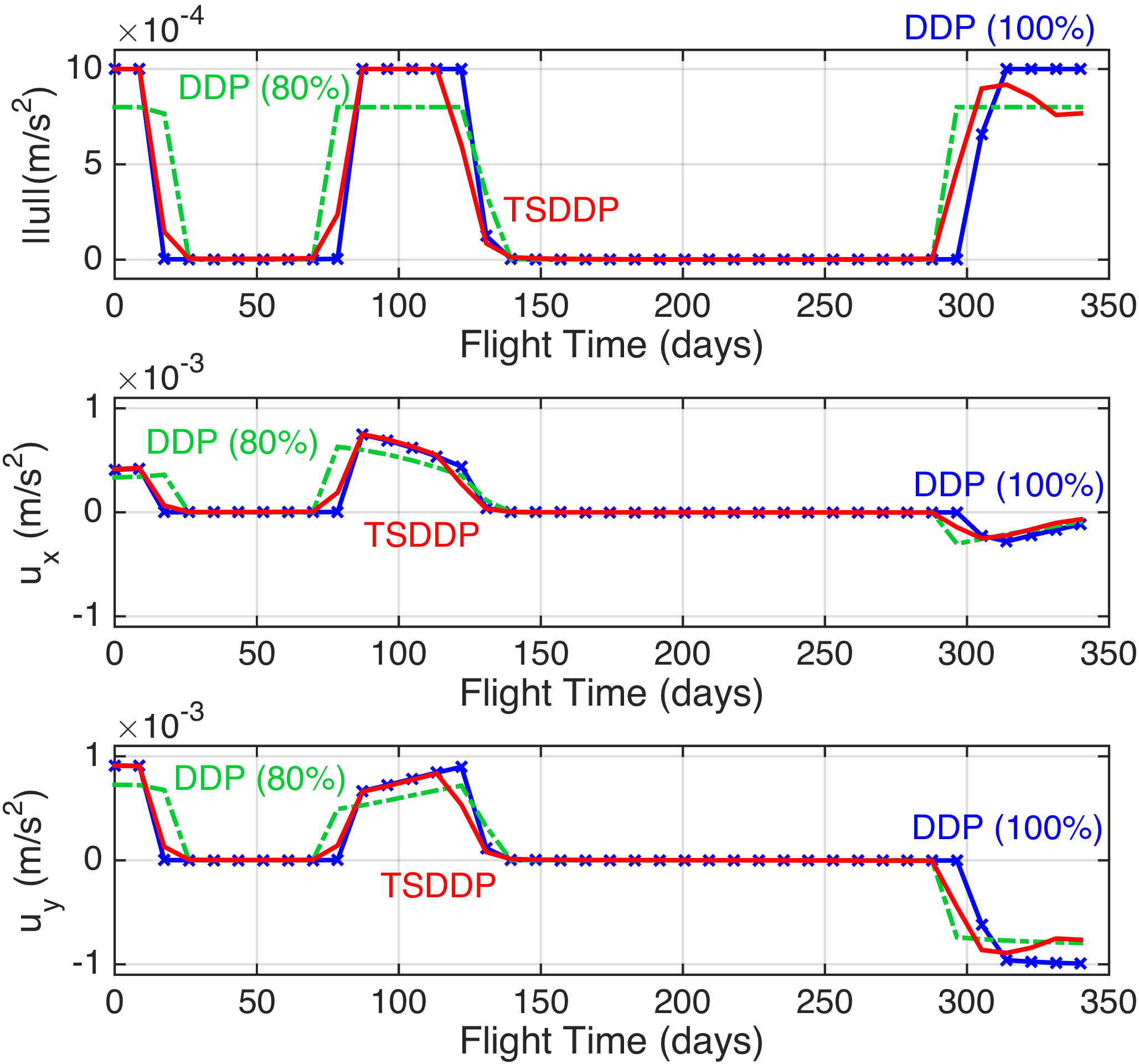}
			{\small b) Nominal thrust profiles}
		\end{center}
      \end{minipage}
    \end{tabular}
    \caption{Nominal trajectory and control inputs.}
    \label{fig:nom_trjcp_all_low_thrust}
  \end{center}
\end{figure}

Figure \ref{fig:nom_trjcp_all_low_thrust} shows the nominal trajectories and control computed by DDP and TSDDP. Comparing DDP with 100\% duty cycle and TSDDP reveals that TSDDP shifts the second thrusting arc ahead, imposes the duty cycle and changes the thrust direction at the third thrusting arc. This result implies that the proposed approach gets a robust solution by adding not only the duty cycle but also various forms of margins, such as thrust direction changes and thrusting time shifts.

\begin{figure}[htbp]
  \begin{center}
    \begin{tabular}{c}
	\hspace*{-0.075\hsize}
      \begin{minipage}{0.55\hsize}
		\begin{center}
			\includegraphics[clip,width=\hsize]{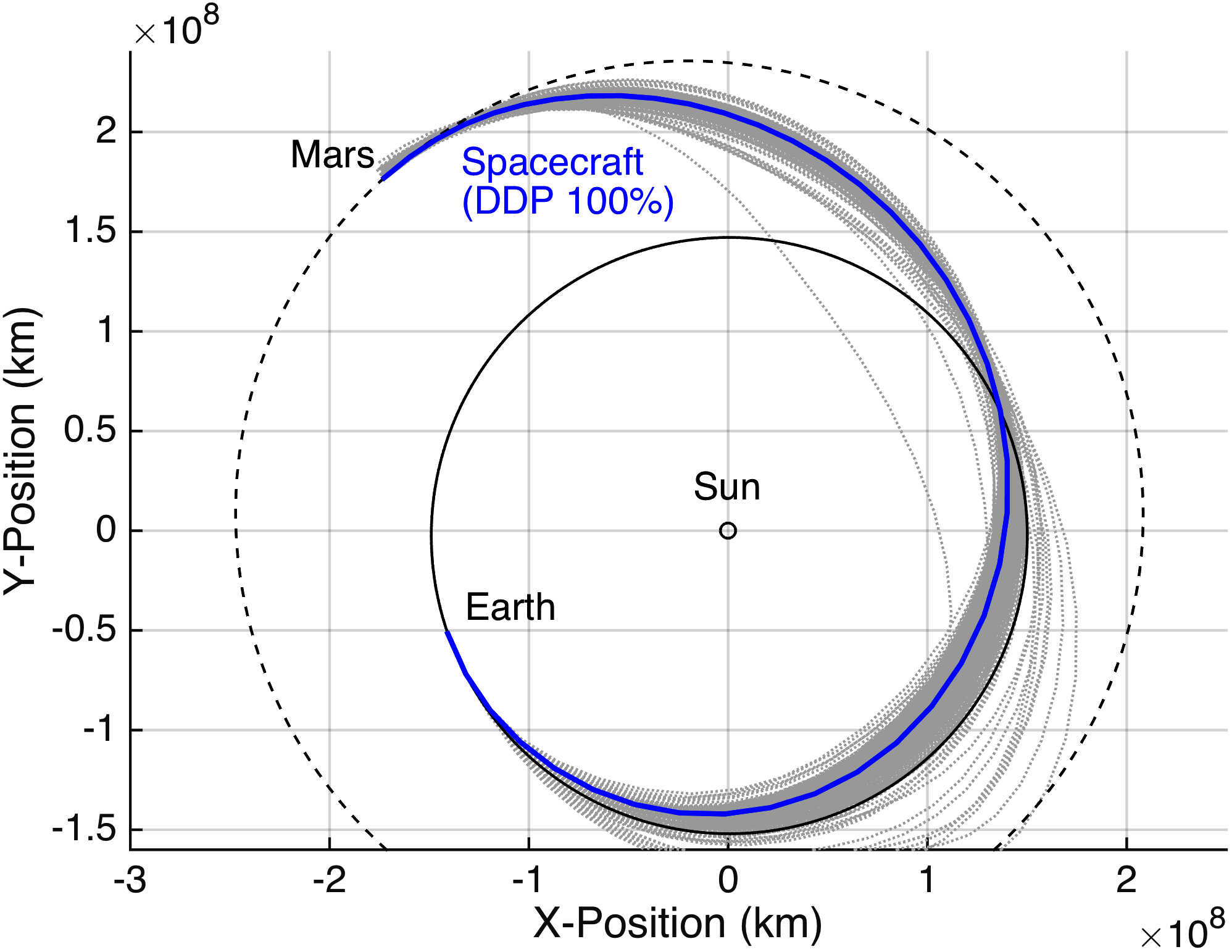}
			{\small a) Sample trajectories (In each sample trajectories, differences from nominal are exaggerated 5 times for illustration purposes.)}
		\end{center}
      \end{minipage}

      \begin{minipage}{0.55\hsize}
		\begin{center}
			\includegraphics[clip,width=\hsize]{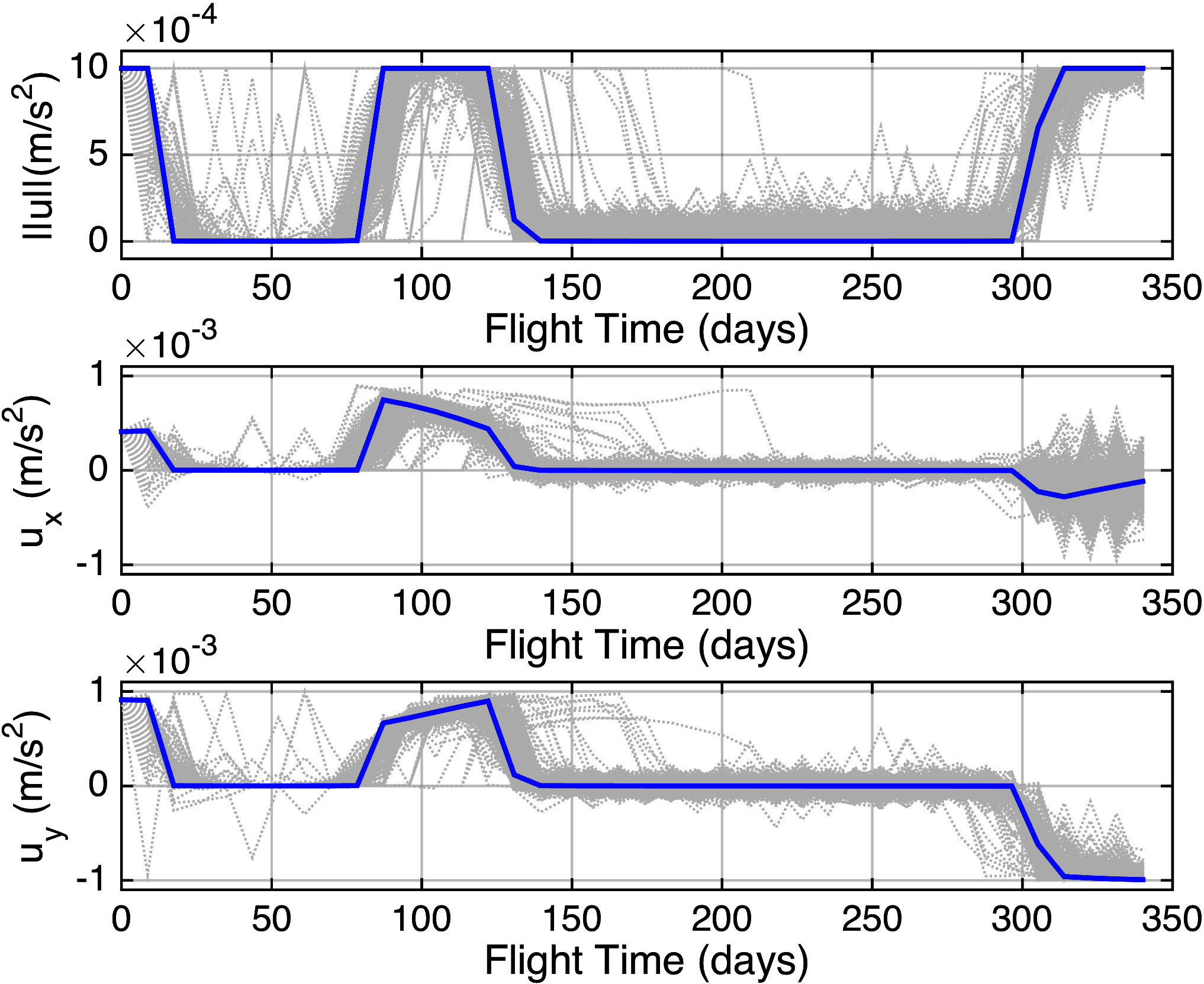}
			{\small b) Sample control profiles}
		\end{center}
      \end{minipage}
    \end{tabular}
    \caption{Sample trajectories and control profiles of DDP with 100\% duty cycle: case 1. (500 samples, each gray line indicates a sample, solid line represents the nominal solution.)}
    \label{fig:mc_trjcp_ddp100_low_thrust}
  \end{center}
\end{figure}

\begin{figure}[htbp]
  \begin{center}
    \begin{tabular}{c}
	\hspace*{-0.075\hsize}
      \begin{minipage}{0.55\hsize}
		\begin{center}
			\includegraphics[clip,width=\hsize]{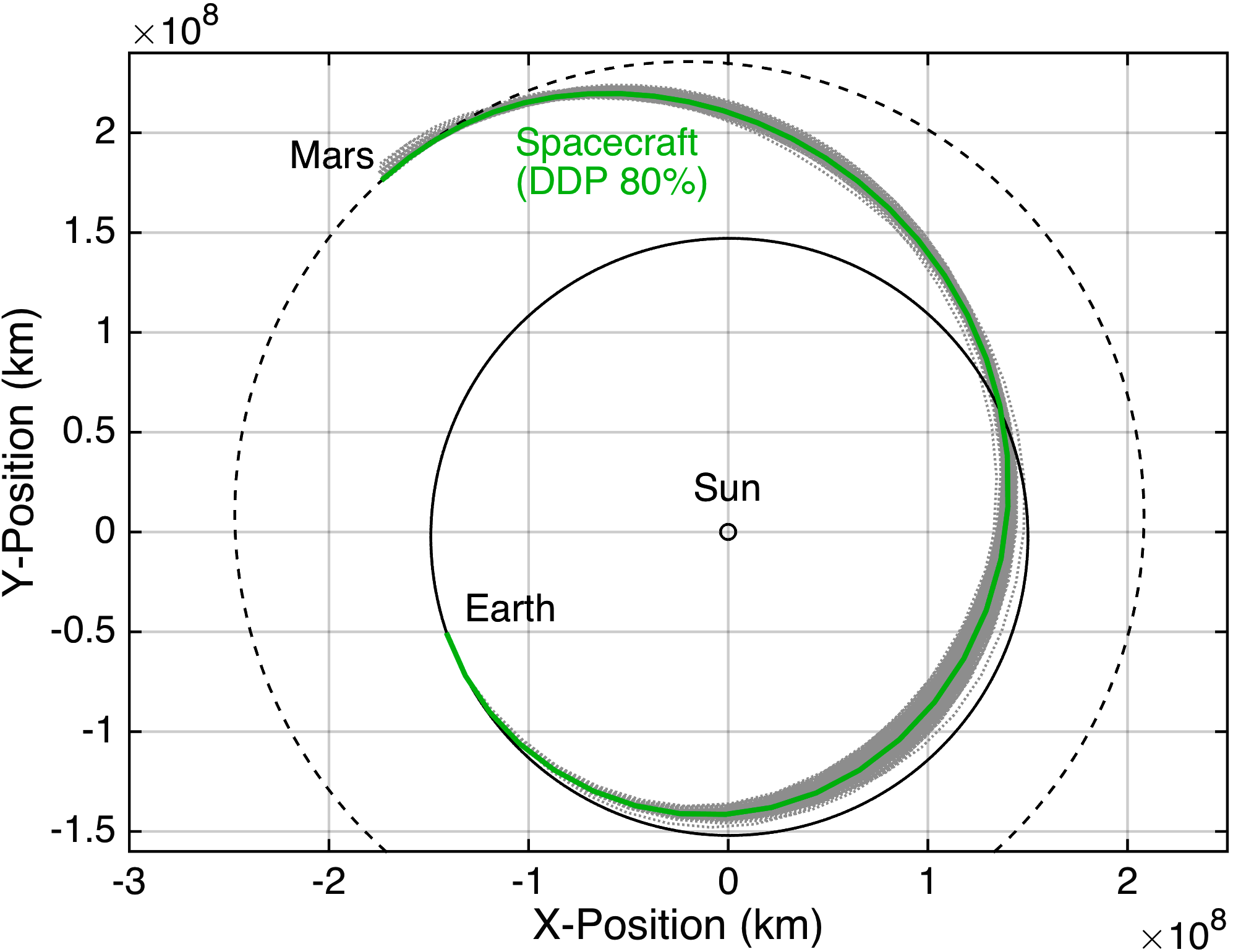}
			{\small a) Sample trajectories (In each sample trajectories, differences from nominal are exaggerated 5 times for illustration purposes.)}
		\end{center}
      \end{minipage}

      \begin{minipage}{0.55\hsize}
		\begin{center}
			\includegraphics[clip,width=\hsize]{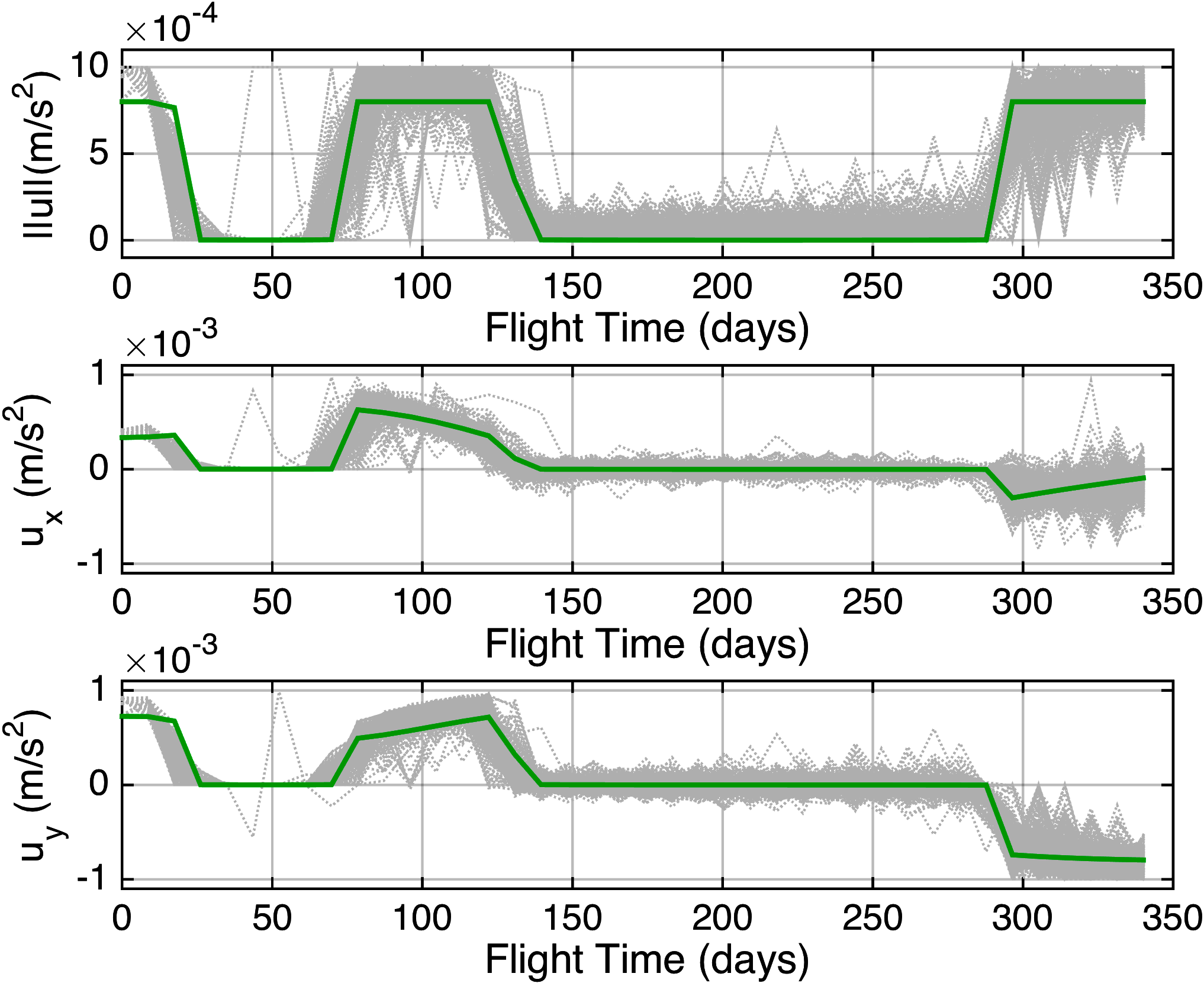}
			{\small b) Sample control profiles}
		\end{center}
      \end{minipage}
    \end{tabular}
    \caption{Sample trajectories and control profiles of DDP with 80\% duty cycle: case 2. (500 samples, each gray line indicates a sample, solid line represents the nominal solution.)}
    \label{fig:mc_trjcp_ddp80_low_thrust}
  \end{center}
\end{figure}

\begin{figure}[htbp]
  \begin{center}
    \begin{tabular}{c}
	\hspace*{-0.075\hsize}
      \begin{minipage}{0.55\hsize}
		\begin{center}
			\includegraphics[clip,width=\hsize]{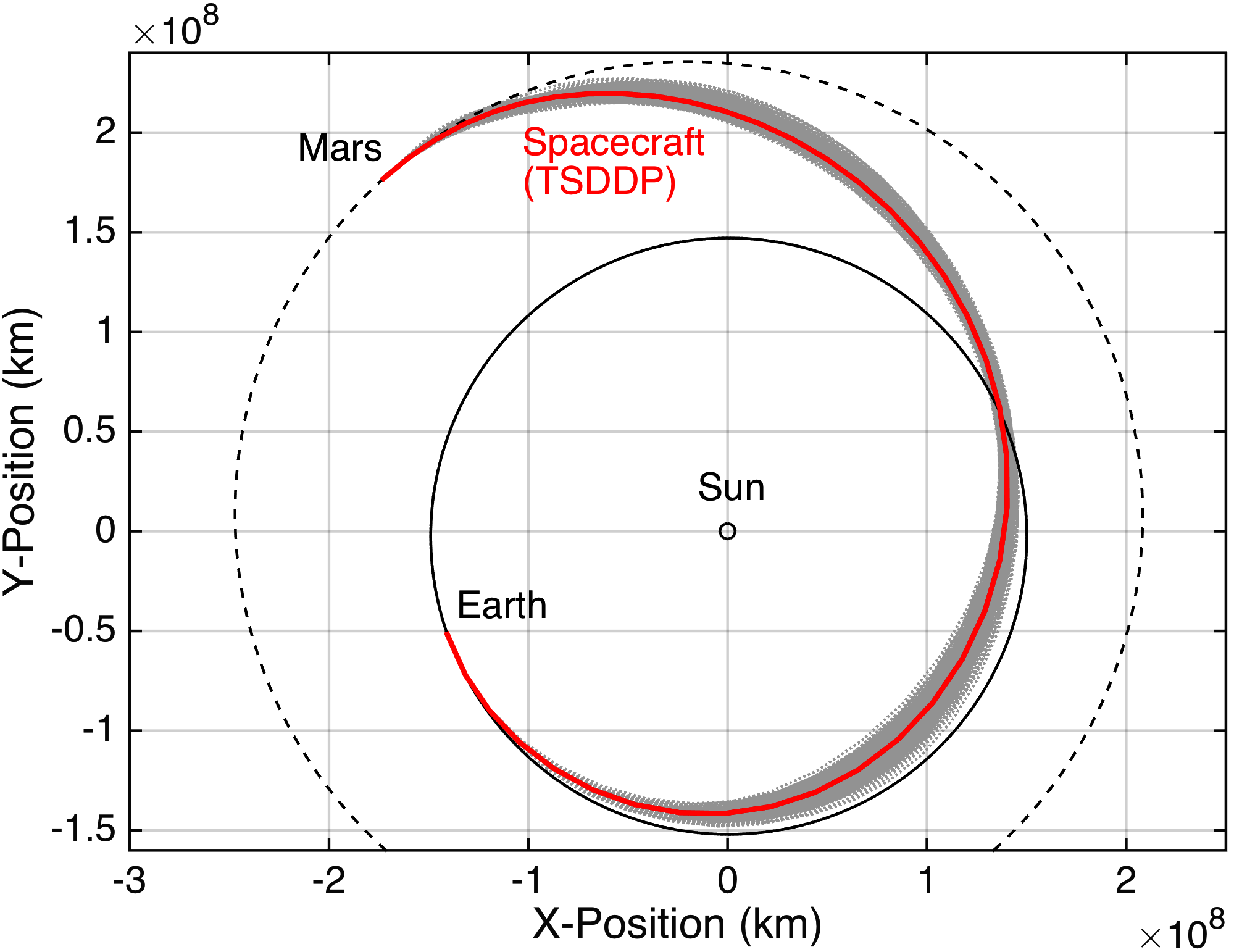}
			{\small a) Sample trajectories (In each sample trajectories, differences from nominal are exaggerated 5 times for illustration purposes.)}
		\end{center}
      \end{minipage}

      \begin{minipage}{0.55\hsize}
		\begin{center}
			\includegraphics[clip,width=\hsize]{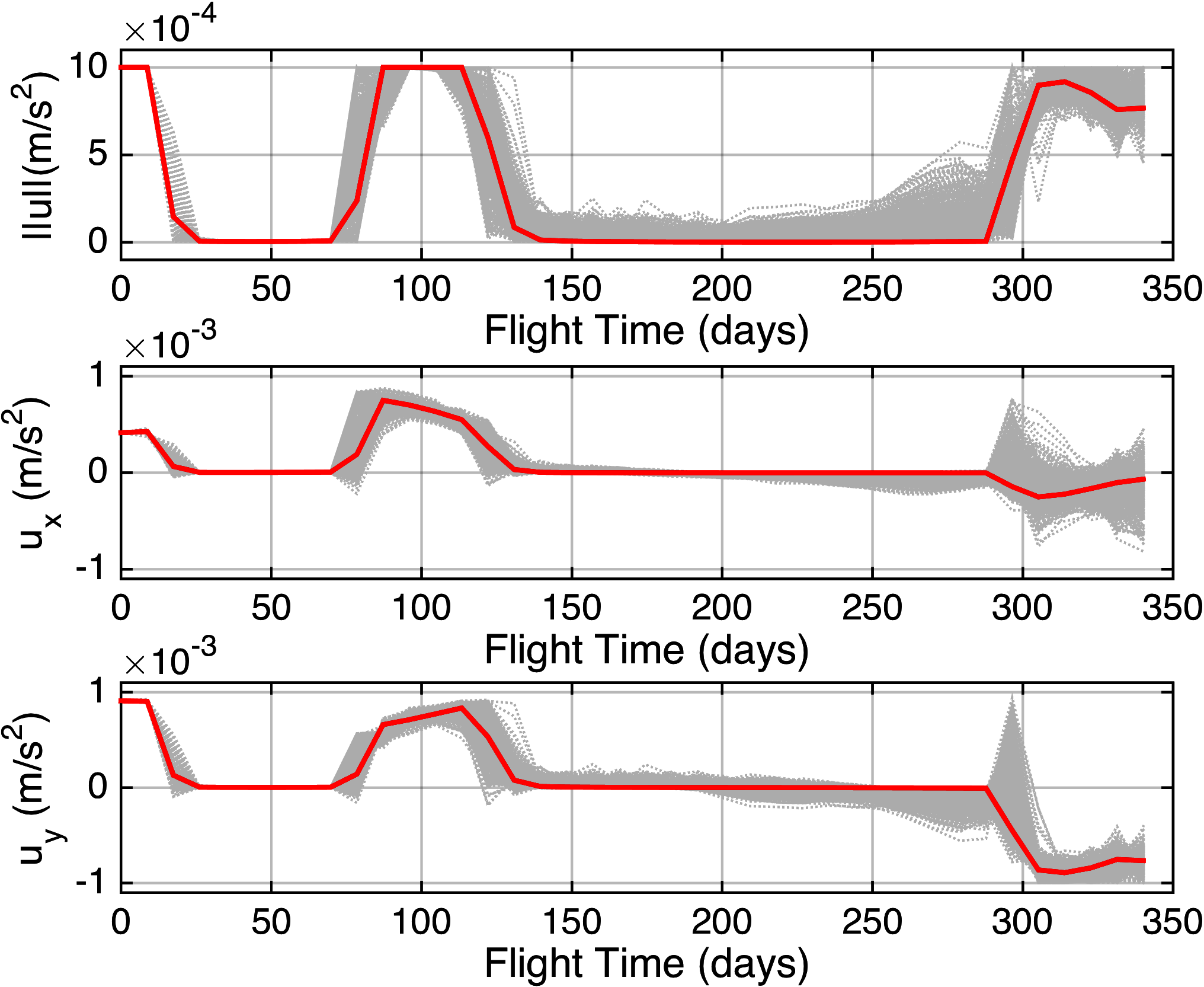}
			{\small b) Sample control profiles}
		\end{center}
      \end{minipage}
    \end{tabular}
    \caption{Sample trajectories and control profiles of TSDDP: case 3. (500 samples, each gray line indicates a sample, solid line represents the nominal solution.)}
    \label{fig:mc_trjcp_tsddp_low_thrust}
  \end{center}
\end{figure}

\begin{figure}[htbp]
	\begin{center}
		\begin{tabular}{c}
			\begin{minipage}{0.55\hsize}
				\begin{center}
					\includegraphics[clip,width=\hsize]{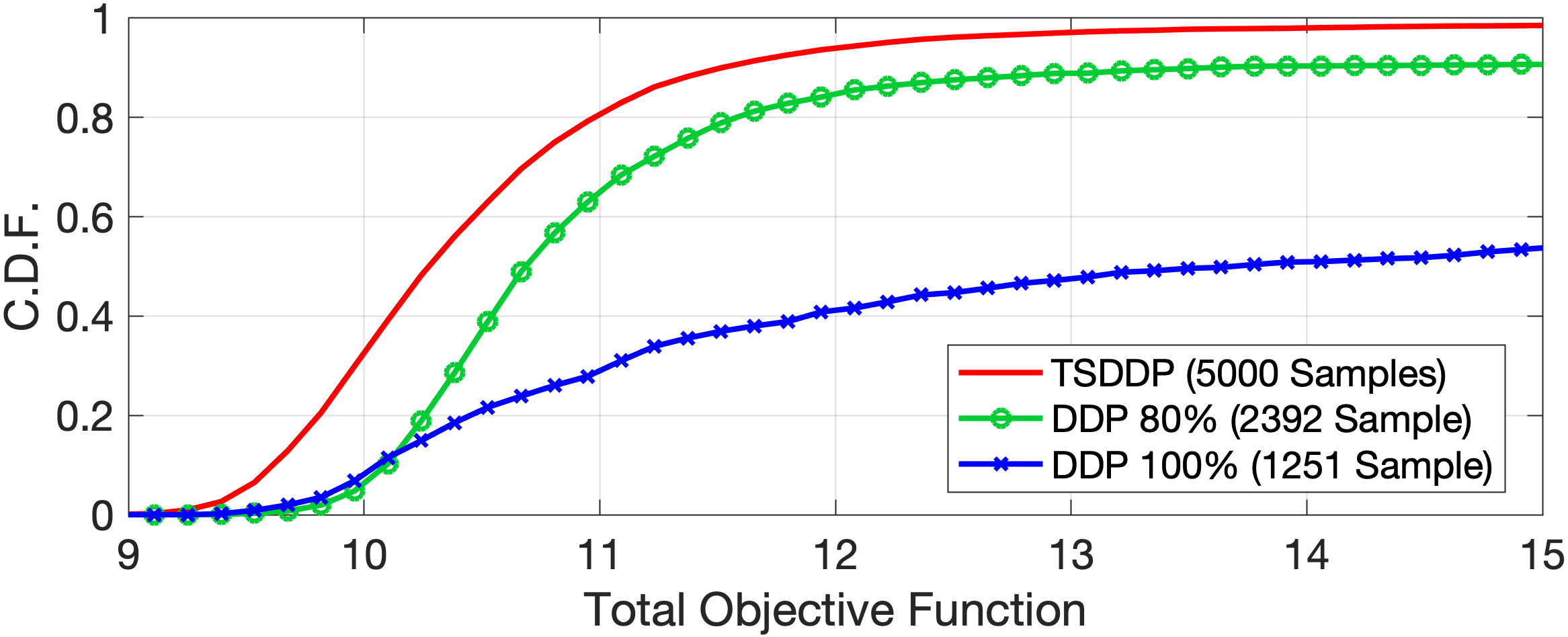}
				\end{center}
			\end{minipage}
			\\
			\begin{minipage}{0.55\hsize}
				\begin{center}
					\includegraphics[clip,width=\hsize]{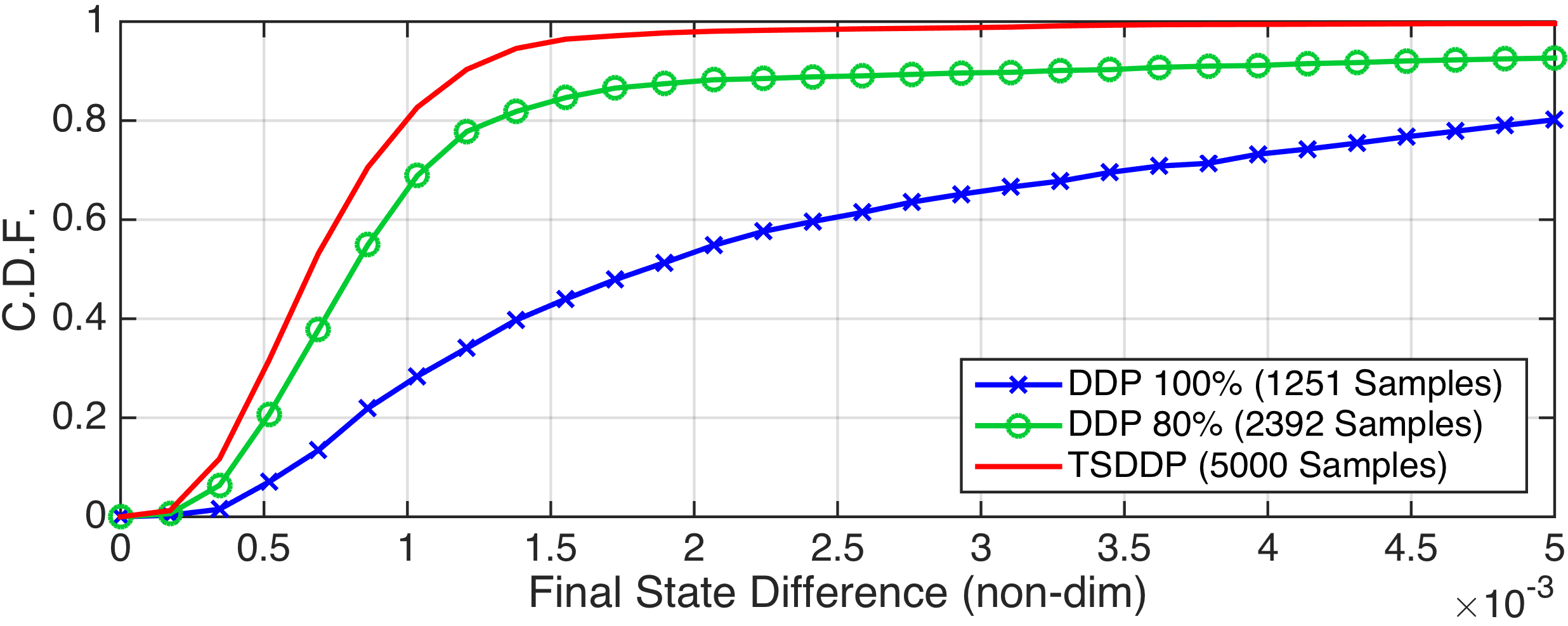}
				\end{center}
			\end{minipage}
			\\
			\begin{minipage}{0.55\hsize}
				\begin{center}
					\includegraphics[clip,width=\hsize]{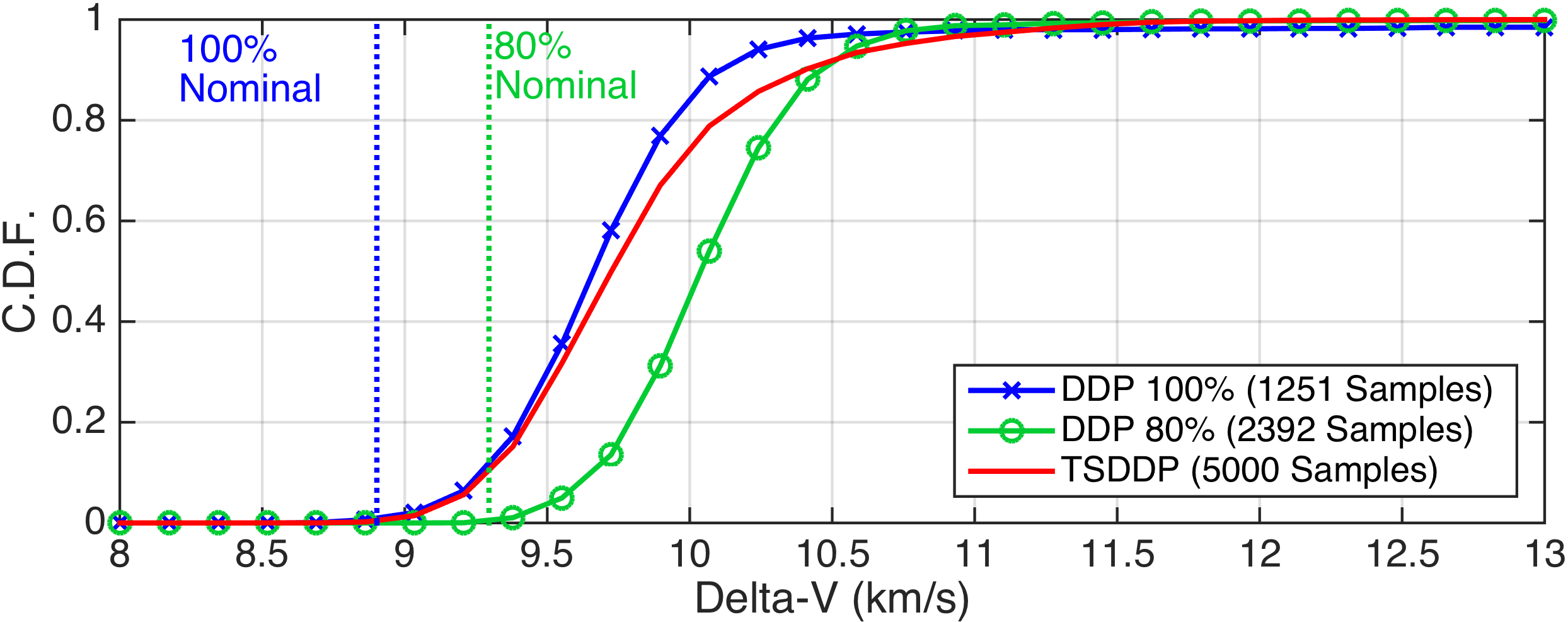}
				\end{center}
			\end{minipage}
		\end{tabular}
		\caption{Cumulative distribution functions (CDFs) of the objective functions.}
		\label{fig:mc_cdf_81p_low_thrust}
	\end{center}
\end{figure}

As with the previous example, we now run Monte Carlo simulations to evaluate the robustness and optimality of DDP with 100\% duty cycle (case 1), DDP with 80\% duty cycle (case 2), and TSDDP (case 3). Although for a fair comparison we should re-optimize all three cases, case 3 requires intensive computational loads. Therefore, in case 3, we adopt a piecewise linear control policy that interpolates the control vectors at the sigma points. The sample trajectories and control profiles are shown in Figs. \ref{fig:mc_trjcp_ddp100_low_thrust} to \ref{fig:mc_trjcp_tsddp_low_thrust}. 
For all cases, most of the sample trajectories do not thrust between the first and second thrusting arc, and most of them correct the trajectories between the second and third thrusting arc. The sample trajectories of case 1 vary drastically from the nominal one, while the trajectories of cases 2 and 3 stay close to the nominal one. One remarkable difference between case 2 and 3 is that case 3 thrusts to the maximum level at the middle of the second thrusting arc, where the energy can be increased most efficiently. Due to this trend, we can expect that case 3 attains a solution with lower $\Delta V$ than case 2.
The robustness and optimality of the solution are examined by CDFs as shown in Fig. \ref{fig:mc_cdf_81p_low_thrust}. As we concluded in the previous example, the two set of the CDFs again show that the proposed method is more robust and optimal than the conventional methods. 

\begin{figure}[htbp]
  \begin{center}
    \begin{tabular}{c}
      \begin{minipage}{0.33\hsize}
		\begin{center}
			\includegraphics[clip,width=\hsize]{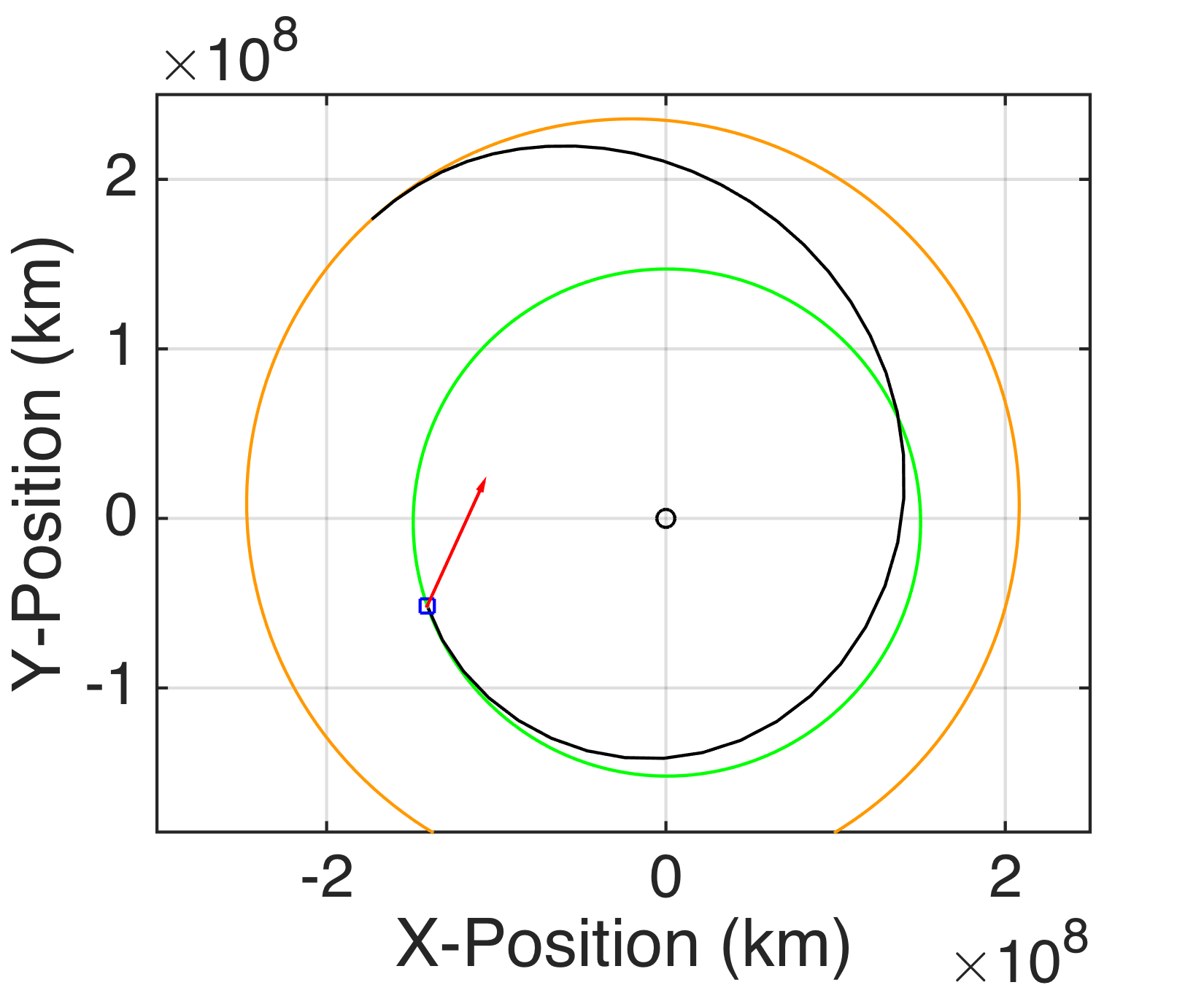}
			{\small a) k=1}
		\end{center}
      \end{minipage}
      \begin{minipage}{0.33\hsize}
		\begin{center}
			\includegraphics[clip,width=\hsize]{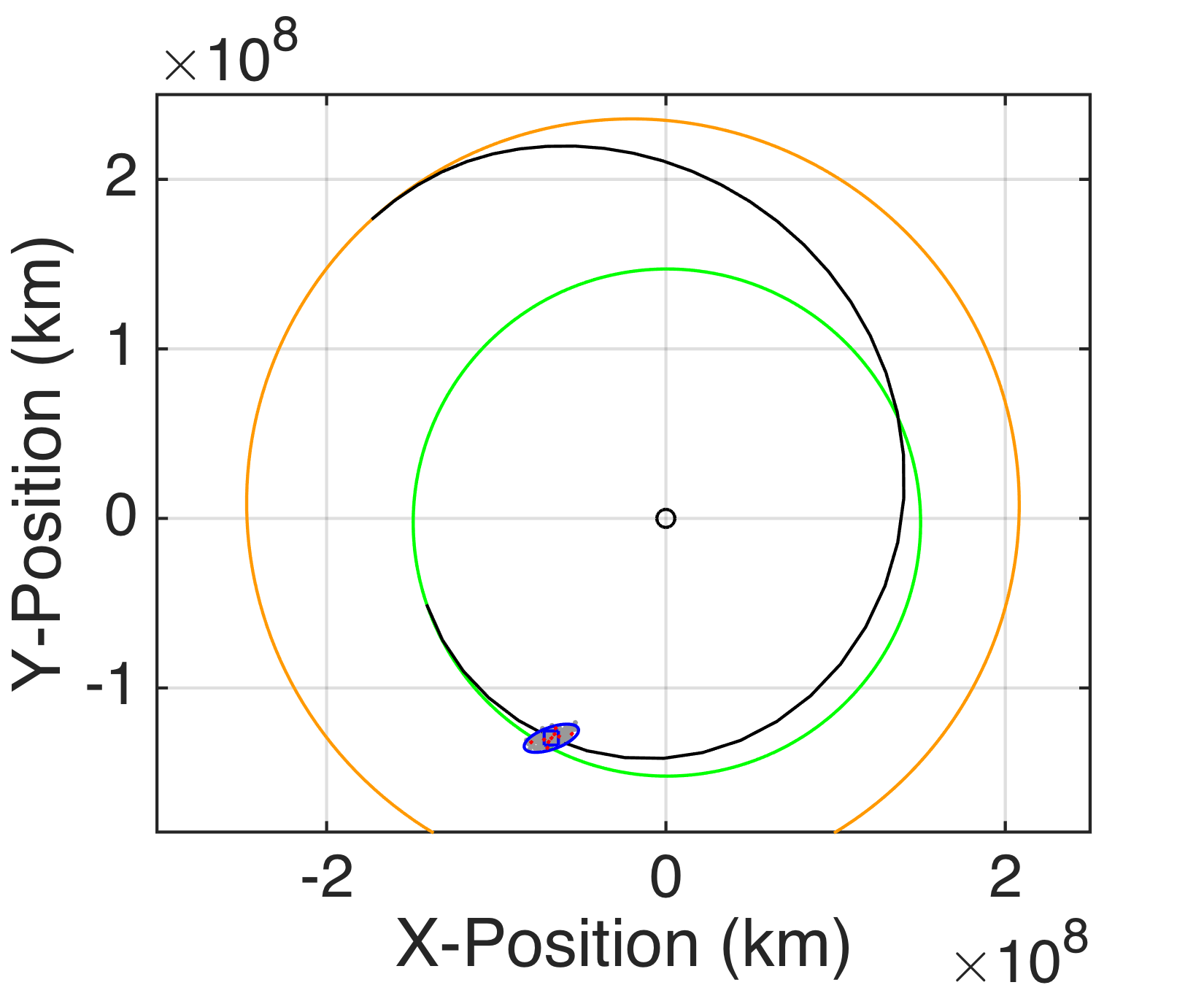}
			{\small b) k=6}
		\end{center}
      \end{minipage}
      \begin{minipage}{0.33\hsize}
		\begin{center}
			\includegraphics[clip,width=\hsize]{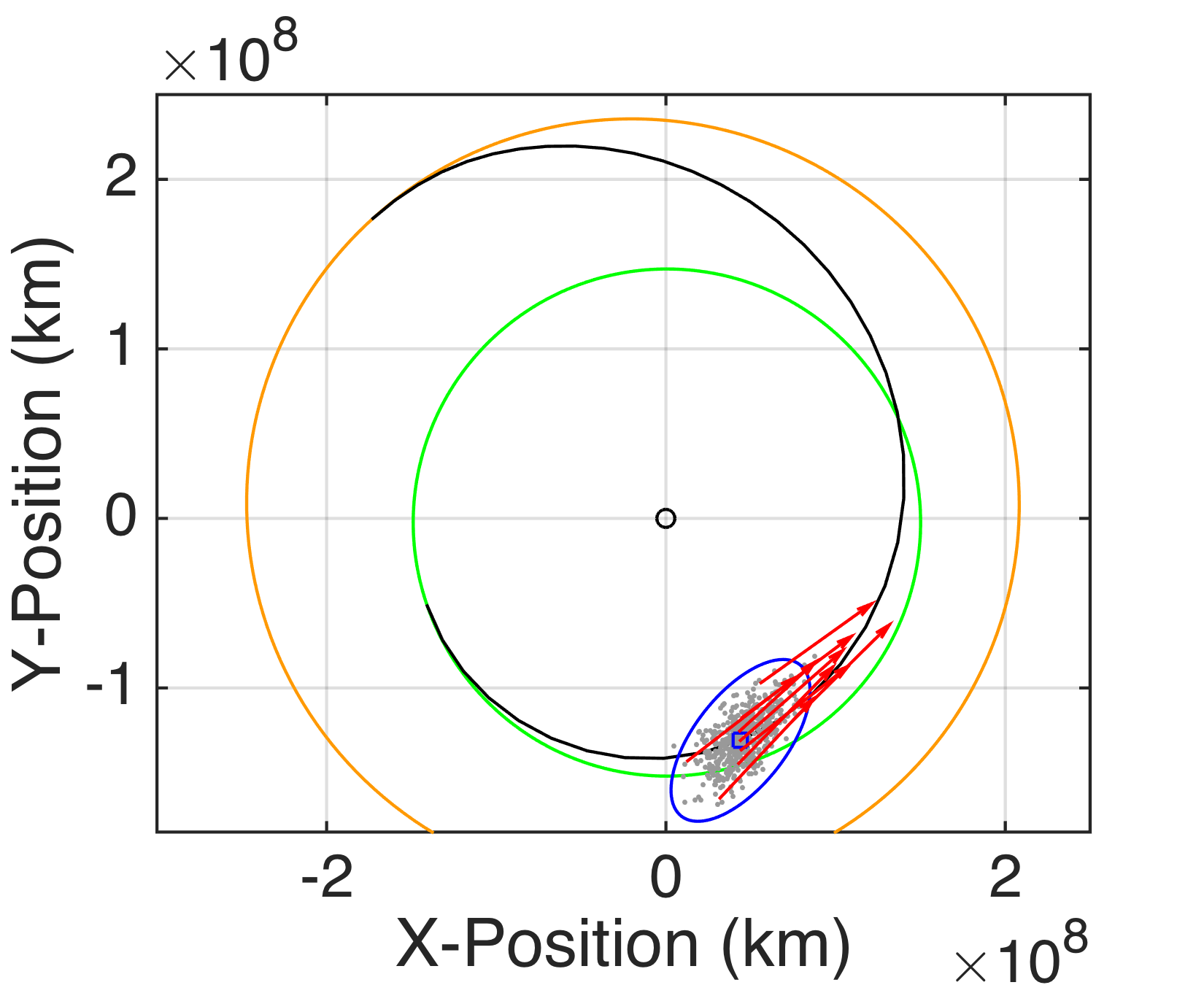}
			{\small c) k=11}
		\end{center}
      \end{minipage}
      \\
      \begin{minipage}{0.33\hsize}
		\begin{center}
			\includegraphics[clip,width=\hsize]{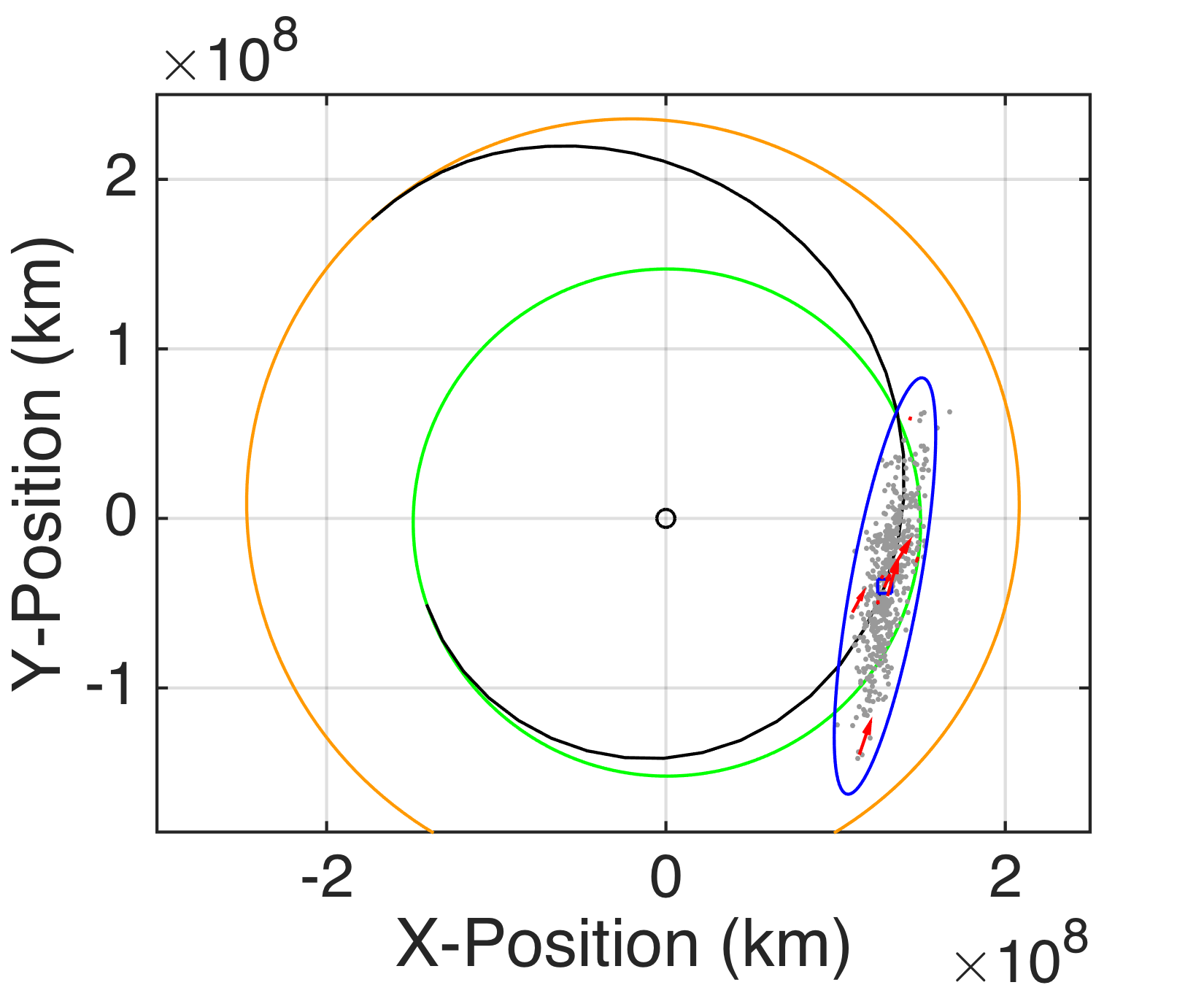}
			{\small d) k=16}
		\end{center}
      \end{minipage}
      \begin{minipage}{0.33\hsize}
		\begin{center}
			\includegraphics[clip,width=\hsize]{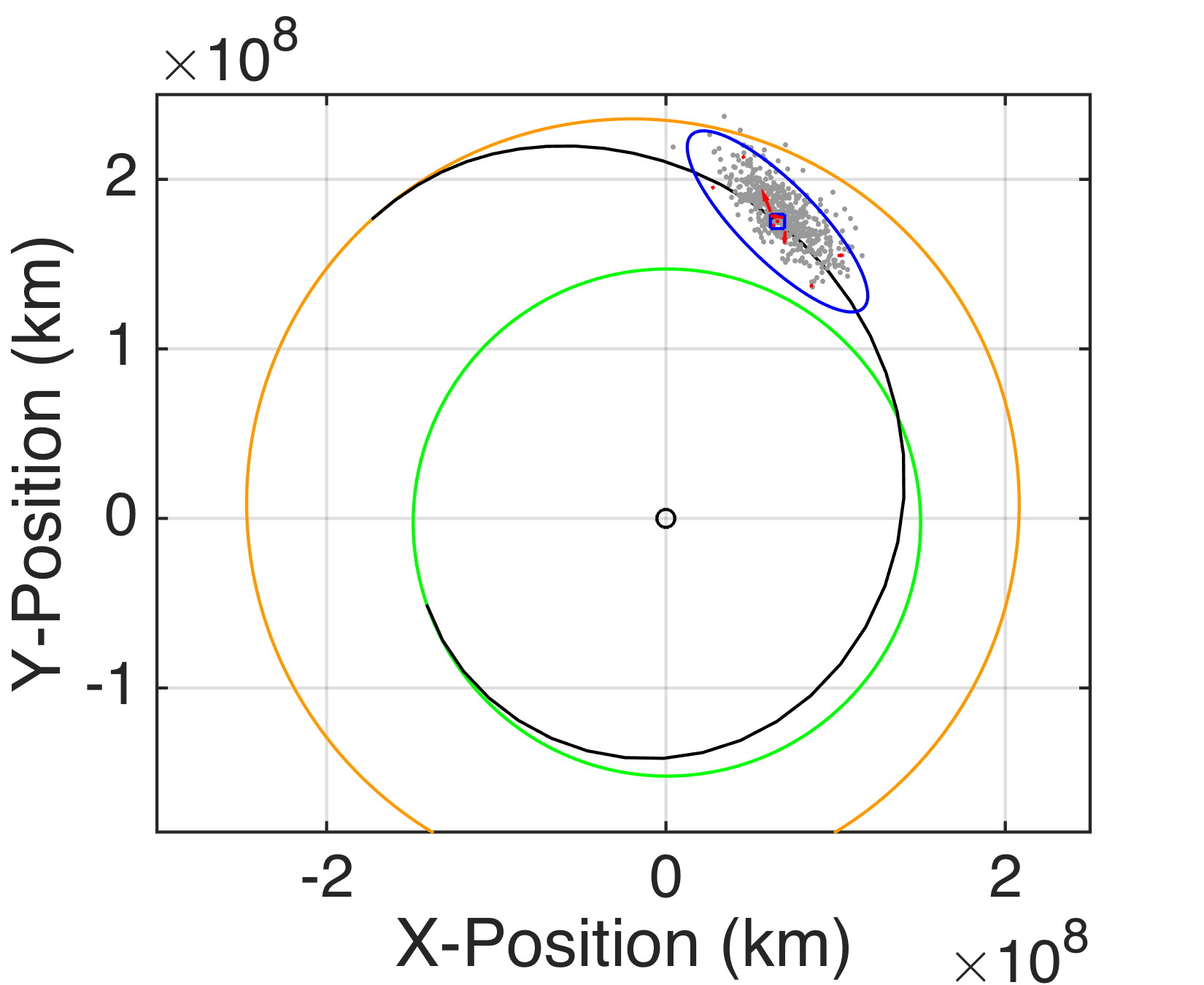}
			{\small e) k=26}
		\end{center}
      \end{minipage}
      \begin{minipage}{0.33\hsize}
		\begin{center}
			\includegraphics[clip,width=\hsize]{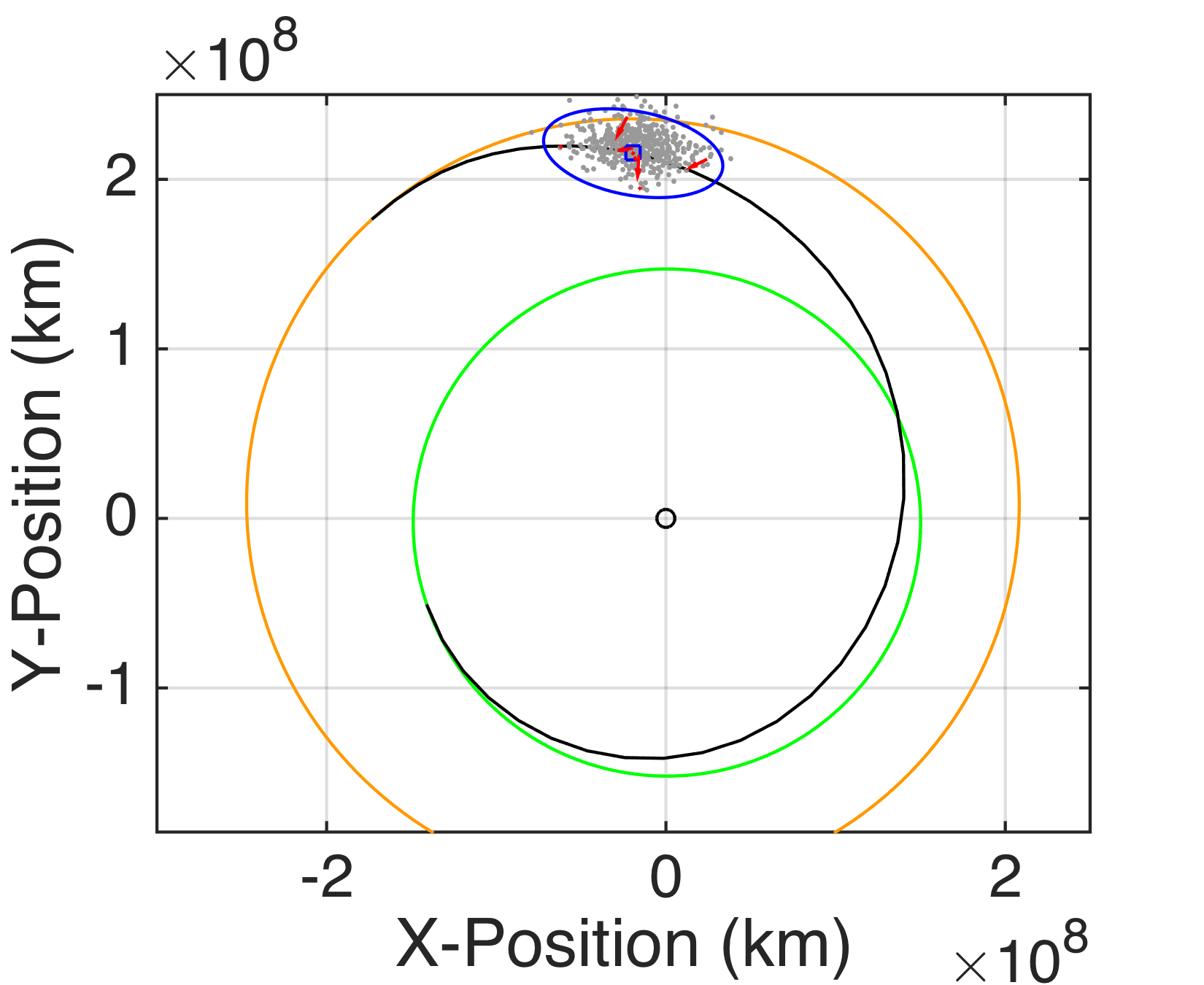}
			{\small f) k=31}
		\end{center}
      \end{minipage}
      \\
      \begin{minipage}{0.33\hsize}
		\begin{center}
			\includegraphics[clip,width=\hsize]{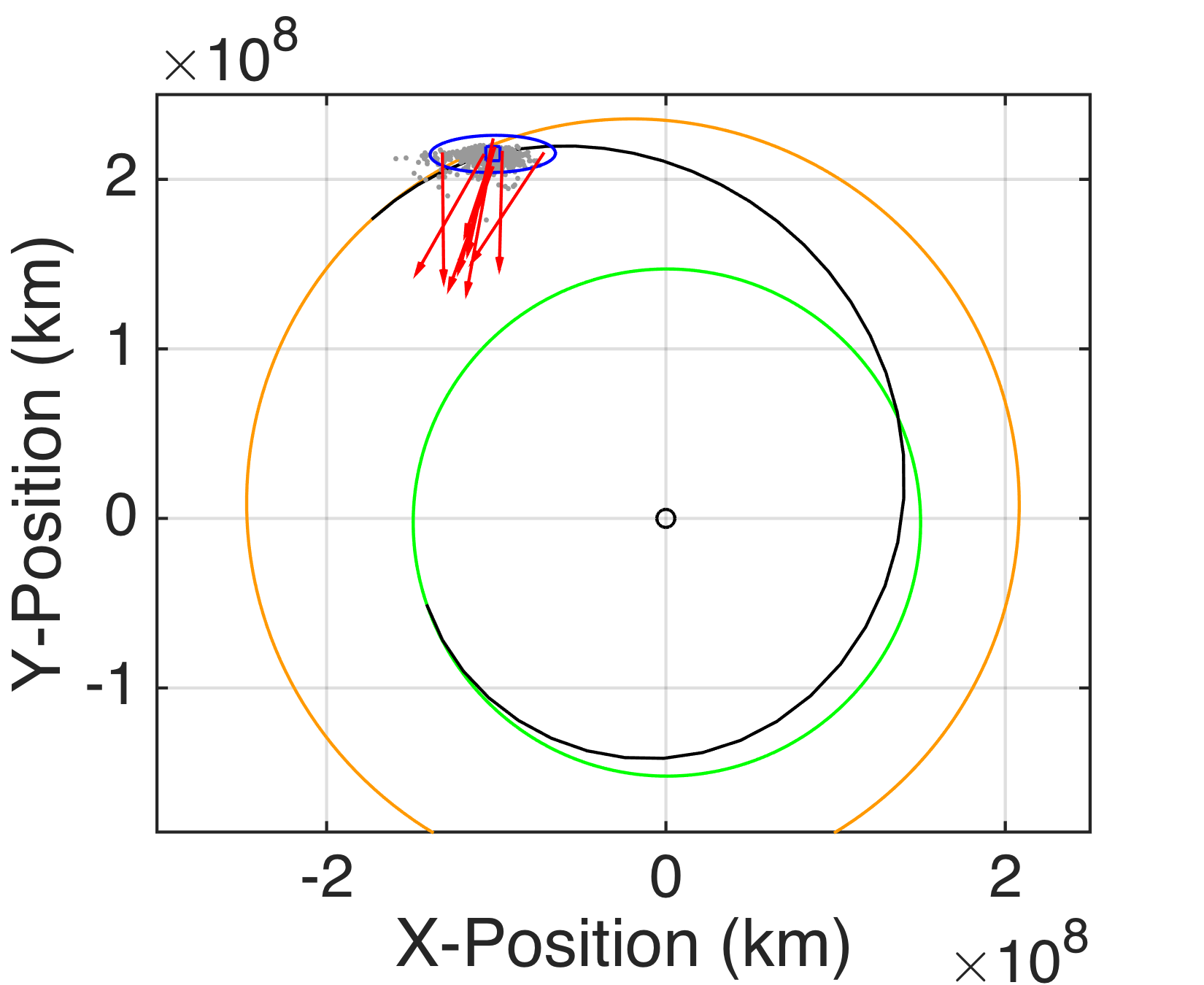}
			{\small g) k=36}
		\end{center}
      \end{minipage}
      \begin{minipage}{0.33\hsize}
		\begin{center}
			\includegraphics[clip,width=\hsize]{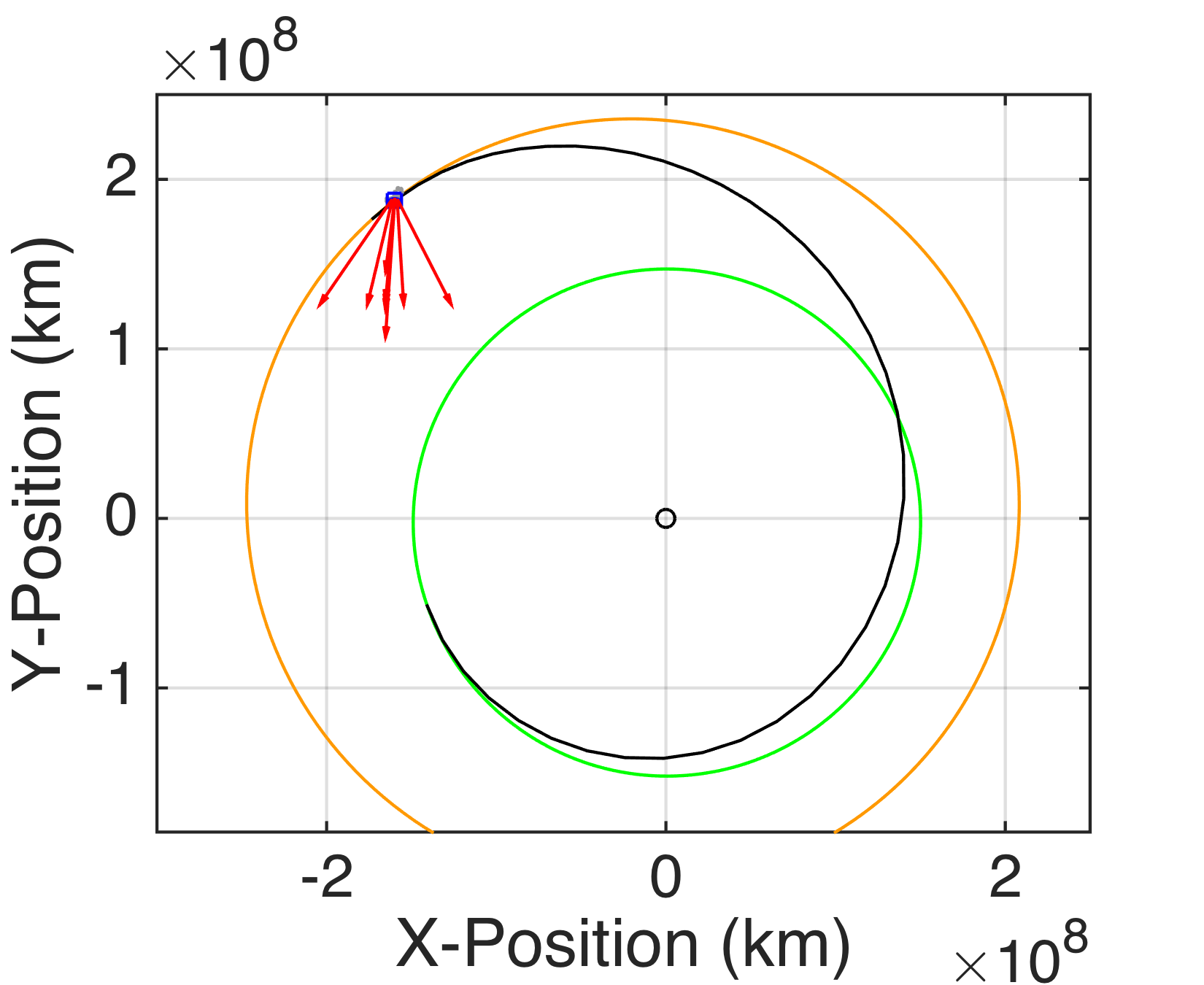}
			{\small h) k=40}
		\end{center}
      \end{minipage}
      \begin{minipage}{0.33\hsize}
		\begin{center}
			\includegraphics[clip,width=\hsize]{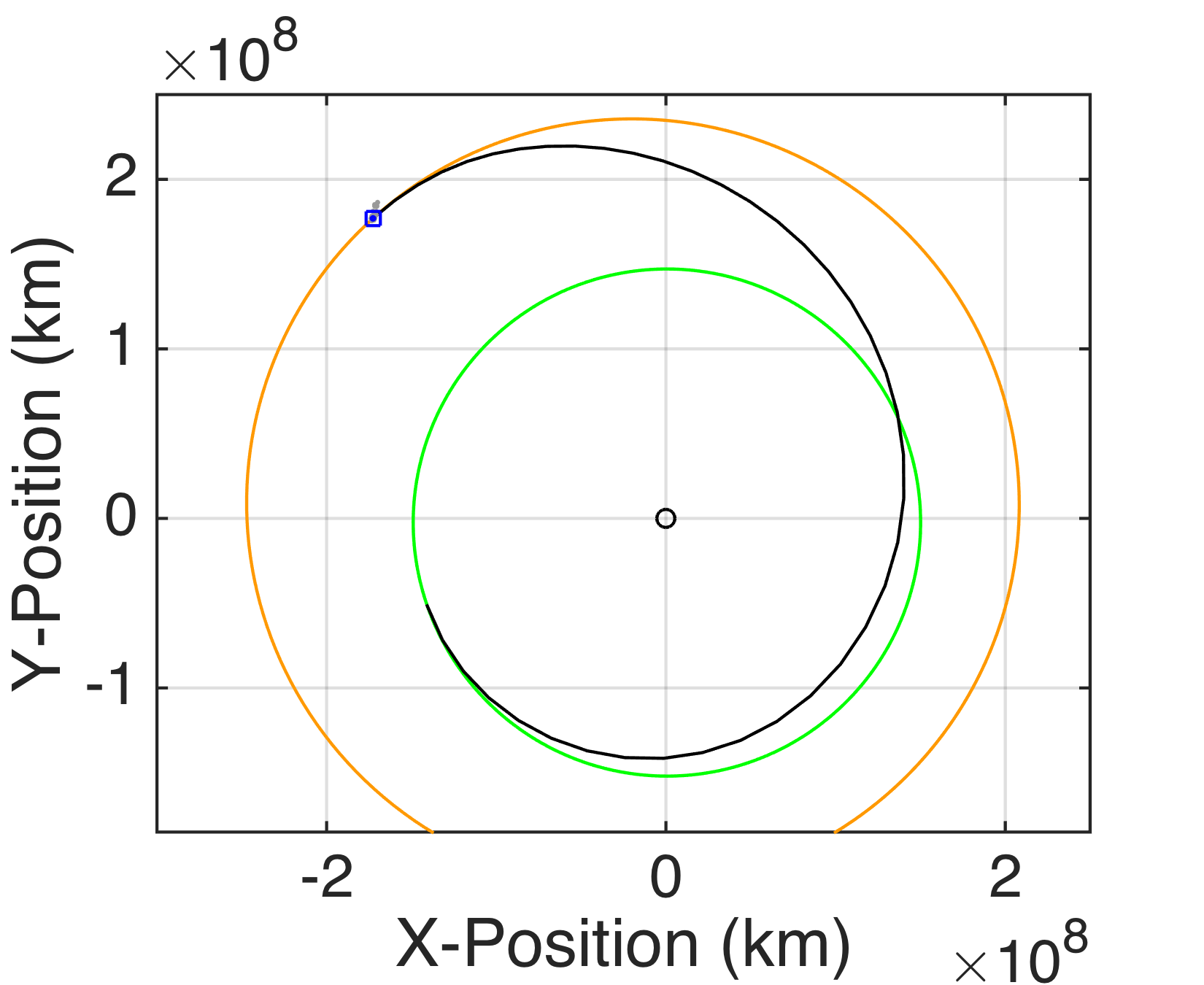}
			{\small i) k=41}
		\end{center}
      \end{minipage}
      
    \end{tabular}
    \caption{Time evolution of TSDDP sample distribution. (Monte Carlo 500 samples, differences from nominal are exaggerated 20 times for illustration purposes, $k$: stage, red arrow: control vectors at sigma points, blue ellipse: 3$\sigma$, gray dots: samples.)}
    \label{fig:mc_capture_tsddp_low_thrust}
  \end{center}
\end{figure}

Figure \ref{fig:mc_capture_tsddp_low_thrust} captures the time evolution of the TSDDP sample distribution.
The Monte Carlo analyses propagate the sample trajectories under a fully nonlinear dynamical system that causes non-Gaussian distributions. Hence, Fig. \ref{fig:mc_capture_tsddp_low_thrust} justifies the approximation of the probability distribution by a Gaussian distribution, and also shows the limitation of the proposed method. For example, in d) $k=16$ in Fig. \ref{fig:mc_capture_tsddp_low_thrust}, most of the samples are enclosed within the 3$\sigma$ ellipse. However, in f) $k=31$ in Fig. \ref{fig:mc_capture_tsddp_low_thrust}, some of them are outside of the 3$\sigma$ ellipse. Future work will introduce high-fidelity uncertainty quantification methods and improve the accuracy of the uncertainty propagation in a non-Gaussian manner.
Note that the simulation in Fig. \ref{fig:mc_capture_tsddp_low_thrust} can be rapidly run once we solve TSDDP and obtain the control policies. Therefore, these results also imply that TSDDP enables robust and optimal autonomous orbital guidance.

%% file: contents/6_conclusion.tex
%
%
\section{Conclusions}

This paper proposes a new method to solve nonlinear constrained stochastic optimal control problems. The method approximates the stochastic process as a Gaussian process and employs the unscented transform to make a stochastic optimal control problem into a deterministic problem. The two numerical examples illustrate the robustness and optimality of the proposed method. The results show that the new method gains robustness by adding various forms of margins such as the duty cycle, thrust direction changes, and thrusting time shifts. Finally, Monte Carlo simulations conclude that the proposed method is more robust and optimal than the conventional methods.

%% file: contents/9_appendix.tex
%
%
\appendix

%
%
\section*{Appendix A: Unscented Transform}\label{sec:app_unscented_transform}
\renewcommand{\theequation}{A\arabic{equation}}
\setcounter{equation}{0}


The unscented transform is a method to calculate the mean value and the covariance of a probability distribution of a random variable that undergoes a nonlinear mapping \cite{Julier1996, Julier1997}. The unscented transform estimates these moments by using a set of representative points, called \textit{sigma points}.

Given the $n$-dimensional random variable $\bm{x}\sim \mathcal{N}(\bar{\bm{x}}, \bm{\mathrm{P}}_x)$, and a nonlinear mapping $\bm{y}=\bm{f}(\bm{x})$, the unscented transform estimates the mean value $\bar{\bm{y}}$ and covariance $\bm{\mathrm{P}}_y$ of the random variable $\bm{y}$ by the following steps.

\begin{enumerate}
	\item Calculate the sigma points $\bm{\mathcal{X}}^{(i)}$ and their weights $c^{(i)}$ for $i\in\mathbb{N}_{2n}$ as
	\begin{align}
	\bm{\mathcal{X}}^{(0)} &= \bar{\bm{x}},\\
	\bm{\mathcal{X}}^{(j)} &= \bar{\bm{x}} + \left(\sqrt{(n+\kappa)\bm{\mathrm{P}}_x}\right)_j,  \ \ \ j\in\mathbb{N}_{1:n},\\
	\bm{\mathcal{X}}^{(n+j)} &= \bar{\bm{x}} - \left(\sqrt{(n+\kappa)\bm{\mathrm{P}}_x}\right)_j, \ \ \ j\in\mathbb{N}_{1:n},\\
	c^{(0)} &= \kappa/n+\kappa,\\
	c^{(j)} &= 1/2(n+\kappa),  \ \ \ j\in\mathbb{N}_{1:2n}
	\end{align}
	where $(\cdot)_j$ represents $j$th column vector of the matrix, the matrix square roots are computed by eigendecomposition, and the arbitrary parameter $\kappa\in(0,\infty)\subset \mathbb{R}$.
	
	\item Obtain the set of the transformed sigma points $\bm{\mathcal{Y}}^{(i)}$
	\begin{equation}
	\bm{\mathcal{Y}}^{(i)} = \bm{f}(\bm{\mathcal{X}}^{(i)}), \ \ \ i\in\mathbb{N}_{2n}
	\end{equation}
	
	\item Calculate the mean value and the covariance of $\bm{y}$ by using the weights $c^{(i)}$ and the transformed sigma points $\bm{\mathcal{Y}}^{(i)}$
	\begin{align}
	\bar{\bm{y}} &= \sum_{i=0}^{2n} c^{(i)} \bm{\mathcal{Y}}^{(i)}\\
	\bm{\mathrm{P}}_y &= \sum_{i=0}^{2n} c^{(i)} \left(\bm{\mathcal{Y}}^{(i)} - \bar{\bm{y}}\right) \left(\bm{\mathcal{Y}}^{(i)} - \bar{\bm{y}}\right)^T
	\end{align}
\end{enumerate}

%
%
\section*{Appendix B: Chance-Constrained Programming}
\renewcommand{\theequation}{B\arabic{equation}}
\setcounter{equation}{0}

Chance-constrained programming is a method to solve stochastic optimization with constraints.\cite{Charnes1959, Li2008, Ono2015} The general formulation of an optimization problem under uncertainty is
\begin{equation}
\min_{\bm{x}} J(\bm{x}, \bm{w})
\end{equation}
subject to
\begin{align}
\bm{c}_{eq}(\bm{x}, \bm{w}) &= \bm{0}\\
\bm{c}(\bm{x}, \bm{w}) &\leq \bm{0}\label{eq:chance_constrained_ineq}
\end{align}
where $J(\cdot)$ is the objective function, $\bm{c}_{eq}(\cdot)$ are the equality constraints, and $\bm{c}(\cdot)$ are the inequality constraints. Also, $\bm{x}$ is the decision variable and $\bm{w}$ is the uncertainty vector. 

The inequality [Eq.(\ref{eq:chance_constrained_ineq})] may not be satisfied due to the uncertainty vector $\bm{w}$, thus the chance-constrained method quantifies the possible violations and modifies the inequality constraint as
\begin{equation}
\mathbb{P}[ \bm{c}(\bm{x}, \bm{w})\leq \bm{0} ] \geq p
\end{equation}
where $p\in[0,1]$ is the probability level.